\documentclass[12pt]{article}
\usepackage{amssymb}
\usepackage{amsmath}
\usepackage{amsthm}
\usepackage{color}
\usepackage{comment}
\usepackage{enumitem}
\usepackage{graphicx}
\def\cF{\mathcal{F}}
\def\cI{\mathcal{I}}
\def\cN{\mathcal{N}}
\def\cO{\mathcal{O}}
\def\ee{\varepsilon}
\newcommand{\removableFootnote}[1]{}

\begin{document}

\title{
The influence of localised randomness on regular grazing bifurcations with applications to impacting dynamics.
}
\author{
D.J.W.~Simpson$^{\dagger}$ and R.~Kuske$^{\ddagger}$\\\\
$^{\dagger}$Institute of Fundamental Sciences\\
Massey University\\
Palmerston North\\
New Zealand\\\\
$^{\ddagger}$Department of Mathematics\\
University of British Columbia\\
Vancouver, BC\\
Canada
}
\maketitle

\begin{abstract}
This paper concerns stochastic perturbations of piecewise-smooth ODE systems
relevant for vibro-impacting dynamics,
where impact events constitute the primary source of randomness.
Such systems are characterised by the existence of switching manifolds
that divide the phase space into regions where the system is smooth.
The initiation of impacts is captured by a grazing bifurcation,
at which a periodic orbit describing motion without impacts
develops a tangential intersection with a switching manifold.
Oscillatory dynamics near regular grazing bifurcations are described by
piecewise-smooth maps involving a square-root singularity, known as Nordmark maps.
We consider three scenarios where coloured noise only affects impacting dynamics,
and derive three two-dimensional stochastic Nordmark maps
with the noise appearing in different nonlinear or multiplicative ways,
depending on the source of the noise.
Consequently the stochastic dynamics differs between the three noise sources,
and is fundamentally different to that of a Nordmark map with additive noise.
This critical dependence on the nature of the noise
is illustrated with a prototypical one-degree-of-freedom impact oscillator.

\end{abstract}



\section{Introduction}
\label{sec:intro}
\setcounter{equation}{0}

Many vibrating mechanical systems experience undesirable impacts that cause wear or sub-optimal performance.
Occasional impacts may be permissible if they result from running components of the system at high speeds for greater efficiency,
and some impacts are unavoidable such as those due to random or rare events.
In these cases it is important to have a clear understanding of the dynamical behaviour that impacts may induce.
Impacting dynamics is often complicated or chaotic because impacts are highly nonlinear phenomena \cite{AwLa03,BlCz99,Br99,Ib09,WiDe00}\removableFootnote{
Other references on impact oscillators include
\cite{AiGu93,AwKu04,AnPl10,Ba58,BaLa07,CaGi06,CoZa95,DaZh05,DeBa07,DeVa00,DiNo08,DyPa12,FePf98,FrNo00,HoHo07}
\cite{InPa08,Iv00,Ko69,KuBa11,LuHa96,LuZh06,OsDi08,PeVa92,PfKu90,Po00,QiSu06,Qu05,Sh85,Sh85b,ShHo83,ShHo83b,Sz13,TuSh88}.
}.

For instance tubes in heat exchangers vibrate at high fluid velocities
and impact against baffles in place to guide the fluid flow.
Simple mathematical models of heat exchangers
reveal that chaotic dynamics may be created at the onset of recurring impacts \cite{PaLi92,DeFr99}.
Rotating cutters spun at high speeds experience repeated contact loss with the material being cut.
The resulting impacts between the cutter and the material may similarly induce chaotic dynamics \cite{Gr88,Wi97,Ba01}.
In contrast, some mechanical systems use impacts to achieve their function.
Atomic force microscopes measure surface topography and the chemical properties of a sample on the nanoscale
by gently hitting the sample with a vibrating cantilever.
In this context it is important to understand the complex impacting dynamics so that
the inverse problem of describing the sample can be performed effectively \cite{DaZh07,RaMe08,MiDa10}.

Impacts can often be modelled accurately by carefully describing the deformations
that components of the system undergo during impacts \cite{St00,St04}.
However, for the purposes of understanding vibro-impacting dynamics,
such a modelling approach is too cumbersome and a low degree-of-freedom ODE model can be more useful.
Despite the low-dimensionality,
such models have been shown to quantitatively match the experimental data of a variety of impacting systems.
Examples include a cam-follower system involving occasional detachments between the cam and the follower \cite{AlDi07,AlDi09},
a pendulum experiencing near-instantaneous impacts with a solid wall \cite{PiVi04},
and compliant impacts of a steel block with an elastic beam \cite{InPa06,InPa08b,InPa10}\removableFootnote{
Also \cite{BaIn09,PaIn10}.
}.

\begin{figure}[b!]
\begin{center}
\setlength{\unitlength}{1cm}
\begin{picture}(15,3.8)
\put(0,0){\includegraphics[height=3.6cm]{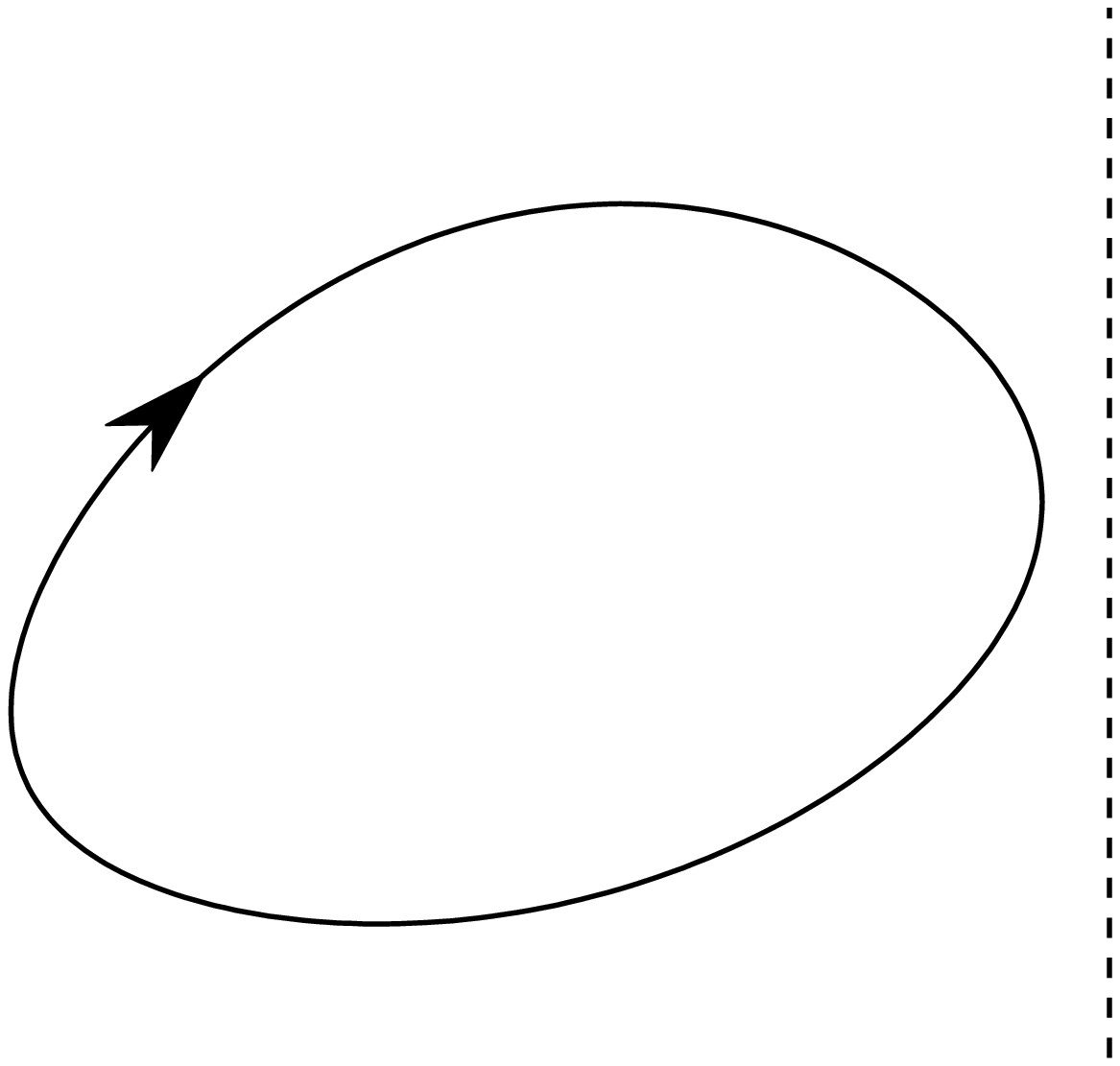}}
\put(5.1,0){\includegraphics[height=3.6cm]{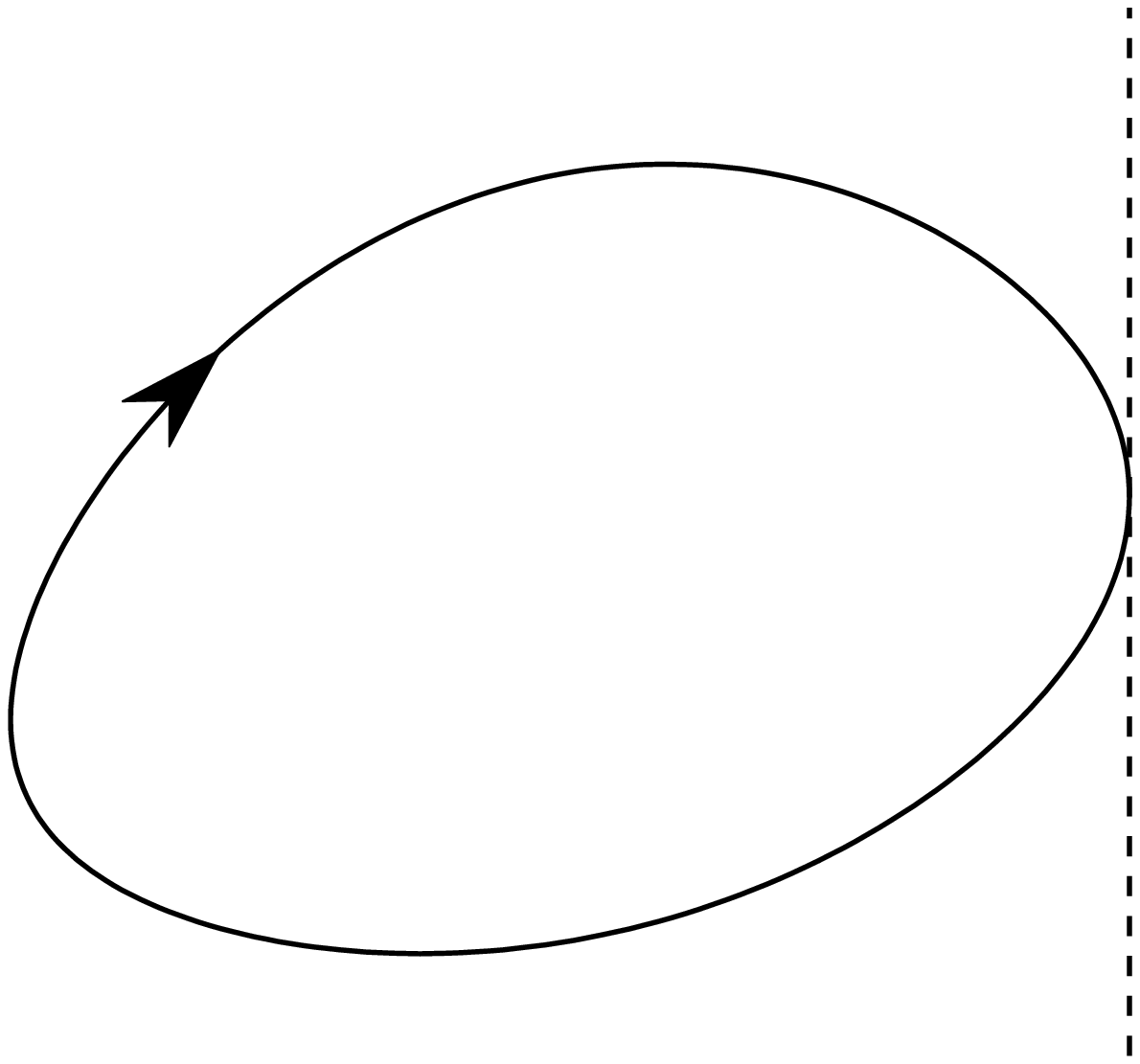}}
\put(10.2,0){\includegraphics[height=3.6cm]{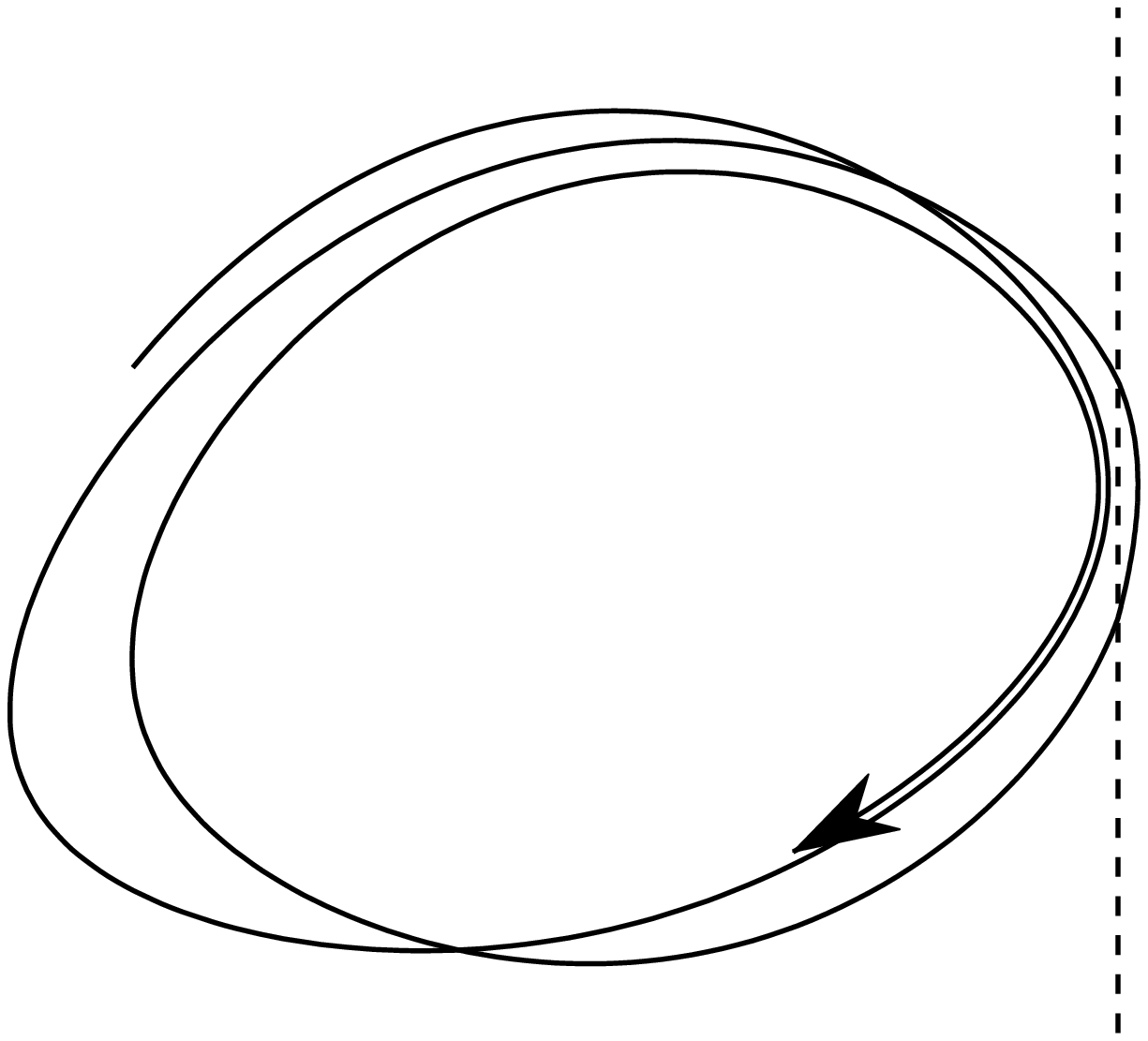}}
\put(.5,3.6){\large \sf \bfseries A}
\put(5.6,3.6){\large \sf \bfseries B}
\put(10.7,3.6){\large \sf \bfseries C}
\put(2,3.6){\small $\eta < 0$}
\put(7.1,3.6){\small $\eta = 0$}
\put(12.2,3.6){\small $\eta > 0$}
\put(4.48,.2){\small $\Sigma$}
\put(9.58,.2){\small $\Sigma$}
\put(14.68,.2){\small $\Sigma$}
\end{picture}
\caption{
Sketches of phase space illustrating a regular grazing bifurcation
occurring at $\eta = 0$, where $\eta \in \mathbb{R}$ is a system parameter.
\label{fig:grazBifSchem}
}
\end{center}
\end{figure}

Using a low degree-of-freedom ODE model, the evolution of the system 
between impacts is tracked in phase space. 
Periodic behaviour without impacts corresponds to a periodic orbit
in phase space that does not reach the switching manifold, $\Sigma$, Fig.~\ref{fig:grazBifSchem}-A.
Physically, $\Sigma$ corresponds to locations or instances where mechanical components come into contact, or lose contact.
As parameters vary, the system may transition from an impact-free regime
to the repeated (though not necessarily regular) occurrence of impacts.
In phase space, the transition occurs when the periodic orbit of the ODE model attains an intersection with $\Sigma$,
Fig.~\ref{fig:grazBifSchem}-B.
This is known as a {\em grazing bifurcation}.

In this paper we study grazing bifurcations of the three-dimensional piecewise-smooth ODE system, 
\begin{equation}
\begin{bmatrix} \dot{u} \\ \dot{v} \\ \dot{w} \end{bmatrix} =
\begin{cases}
f_L(u,v,w;\eta) \;, & u < 0 \\
f_R(u,v,w;\eta) \;, & u > 0
\end{cases} \;,
\label{eq:f}
\end{equation}
where $f_L$ and $f_R$ are smooth functions,
$\eta \in \mathbb{R}$ is a parameter,
and the coordinates $(u,v,w)$ are chosen so that $\Sigma$ is simply the coordinate plane $u=0$.
We assume that for $\eta < 0$, there exists an attracting periodic orbit
describing non-impacting dynamics located entirely in the region $u < 0$,
and that the periodic orbit grazes $\Sigma$ at the origin when $\eta = 0$.
In the context of vibro-impacting systems,
$u < 0$ corresponds to not-in-contact dynamics governed by $f_L$,
and $u > 0$ corresponds to in-contact dynamics governed by $f_R$. 
Impacts may instead be modelled as instantaneous events 
with energy loss and velocity reversal,
in which case a map is usually defined on the switching manifold to describe the action of an impact \cite{VaSc00}.

Theoretical studies of piecewise-smooth and hybrid dynamical systems
have led to a useful classification of grazing bifurcations \cite{DiBu08}.
This paper concerns regular grazing bifurcations.
The grazing bifurcation of (\ref{eq:f}) at $\eta = 0$ is said to be {\em regular} if
\begin{equation}
{\rm sgn} \left( e_1^{\sf T} f_L(0,v,w;\eta) \right) = 
{\rm sgn} \left( e_1^{\sf T} f_R(0,v,w;\eta) \right) \;,
\label{eq:regularGrazingCondition}
\end{equation}
for all $(v,w;\eta)$ in a neighbourhood of $(0,0;0)$.
This condition arises naturally in mechanical systems with compliant impacts,
and implies that $\Sigma$ is neither attracting nor repelling at any point.

As indicated in Fig.~\ref{fig:grazBifSchem}-C,
the steady-state dynamics of (\ref{eq:f}) for $\eta > 0$ is often complicated.
For this reason it is valuable to study the oscillatory dynamics using a return map
based on the points on a Poincar\'{e} section.
A normal form for such a map for regular grazing in $\mathbb{R}^3$ is the Nordmark map 
\begin{equation}
\begin{bmatrix} x_{i+1} \\ y_{i+1} \end{bmatrix} =
\begin{cases}
g(x_i,y_i) \;, & x_i \le 0 \\
g(x_i,y_i) - \begin{bmatrix} \chi \sqrt{x_i} \\ 0 \end{bmatrix} \;, & x_i \ge 0
\end{cases} \;,
\label{eq:N}
\end{equation}
where
\begin{equation}
g(x,y) = \begin{bmatrix} \tau & 1 \\ -\delta & 0 \end{bmatrix}
\begin{bmatrix} x \\ y \end{bmatrix} +
\begin{bmatrix} 0 \\ 1 \end{bmatrix} \mu \;,
\end{equation}
with $\tau, \delta \in \mathbb{R}$, and $\chi = \pm 1$,
as determined by the sign of certain coefficients, see \S\ref{sec:Nordmark}.
The coordinates $(x,y)$ represent points on a Poincar\'{e} section of (\ref{eq:f}),
with $\mu \in \mathbb{R}$ the bifurcation parameter and the grazing bifurcation occurring at $\mu = 0$.
The Nordmark map (\ref{eq:N}), applicable also to models with instantaneous impacts \cite{No91,No97}, 
includes only the leading order terms of the return map, so is valid only for dynamics
close to the grazing bifurcation \cite{DiBu08}.

Each iteration of (\ref{eq:N}) corresponds to one oscillation of (\ref{eq:f}) near the grazing periodic orbit.
The utility of (\ref{eq:N}) lies in the fact that the nature of the dynamics
can be identified by the location of the corresponding points in the $(\mu,x)$-plane.
Specifically, (\ref{eq:N}) is formulated so that
if $x_i < 0$ then the oscillation lies entirely in $u < 0$,
and if $x_i > 0$ then the oscillation enters $u > 0$.
The non-impacting, attracting periodic orbit shown in Fig.~\ref{fig:grazBifSchem}-A
corresponds to the fixed point of $g$,
\begin{equation}
\begin{bmatrix} x^L \\ y^L \end{bmatrix} =
\frac{1}{\delta-\tau+1} \begin{bmatrix} 1 \\ 1-\tau \end{bmatrix} \mu \;,
\label{eq:xLyL}
\end{equation}
with $x^L < 0$ and $\mu < 0$. 

\begin{figure}[b!]
\begin{center}
\setlength{\unitlength}{1cm}
\begin{picture}(12,6)
\put(0,0){\includegraphics[height=6cm]{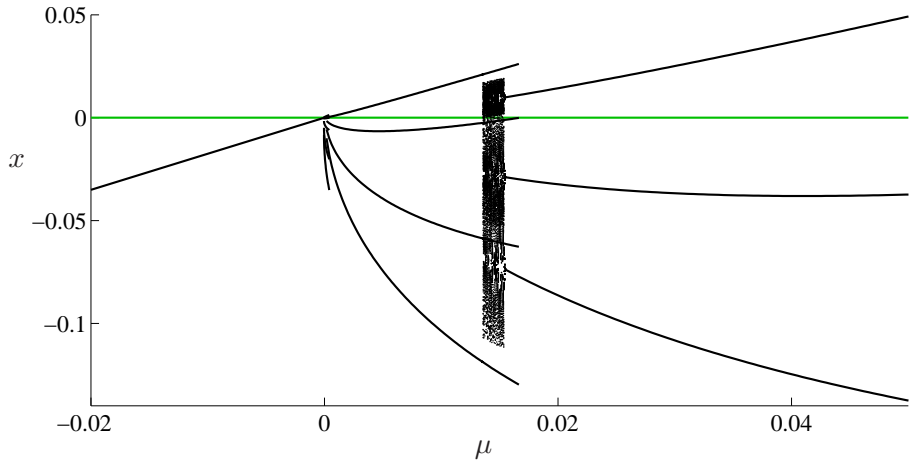}}
\put(6.2,0){\small $\mu$}
\put(0,3.8){\small $x$}
\end{picture}
\caption{
A bifurcation diagram of the Nordmark map (\ref{eq:N}) with
$\tau \approx 0.5813$, $\delta \approx 0.1518$, and $\chi = 1$.
These parameter values correspond to the vibro-impacting system described in \S\ref{sec:oscillator},
with $(k_{\rm osc},b_{\rm osc},k_{\rm supp},b_{\rm supp},d) = (4.5,0.3,10,0,0.1)$.
The fixed point of (\ref{eq:N}) for $\mu < 0$, given by (\ref{eq:xLyL}), corresponds
to an attracting, non-impacting periodic orbit of period $2 \pi$.
The map (\ref{eq:N}) has attracting $3$ and $4$-cycles for certain values of $\mu > 0$ as shown,
as well as an attracting $5$-cycle in the approximate range, $-0.00004 < \mu < 0.00041$.
These correspond to periodic orbits of period approximately equal to $2 k \pi$, for $k = 3,4,5$,
that experience one impact per period.
There also appears to be a chaotic attractor for the approximate range, $0.0135 < \mu < 0.0155$.
\label{fig:detBifDiag}
}
\end{center}
\end{figure}

More generally, a periodic orbit of (\ref{eq:f}) appears 
as a finite set of points in the $(\mu,x)$ plane, 
as in Fig.~\ref{fig:detBifDiag} which shows a typical bifurcation diagram of (\ref{eq:N}).
This figure shows a period-incrementing cascade, corresponding to several different periodic orbits,
and apparently chaotic dynamics, as indicated by a cloud of points\removableFootnote{
I numerically identified the range of parameter values over which the attractor
exists by computing forward orbits from many different initial points, see {\sc goBasinBifDiag.m}.
However, I have not identified the bifurcations associated with the chaotic attractor 
like I did for the corresponding attractor in \cite{SiHo13}.
}.
For alternate values of $\tau$ and $\delta$,
the fixed point $(x^L,y^L)$ may bifurcate directly to chaos \cite{DiBu08,No01}\removableFootnote{
Other dynamics is described in \cite{ZhDa06,DuDe09}.
Direct transitions to chaos at the grazing bifurcation, or to a complicated periodic motion,
have been observed in numerical simulations of impacting systems in \cite{Lu06} and \cite{LuCh07} respectively.
}.
The square-root term in (\ref{eq:N})
is an artifact of the tangency between the periodic orbit and the switching manifold of (\ref{eq:f}) at the grazing bifurcation,
and is responsible for the distinctive shape of the bifurcation diagram near $\mu = 0$.
For some vibro-impacting systems it is more appropriate for the return map to be piecewise-linear
and either continuous \cite{LeVa02} or discontinuous \cite{BuPi06}.
Such maps predict fundamentally different bifurcation structures to those of (\ref{eq:N}).

In order to properly explain complicated vibro-impacting dynamics,
the effects of randomness and uncertainties needs to be taken into account.
Mechanical systems are subject to background vibrations and other sources of noise.
Experimentally measured parameters involve error, and some physical features are left unmodelled.
For instance, one-degree-of-freedom models do not capture high frequency modes
that are often excited by impacts \cite{OeHi97}.

To quantitatively describe stochastic impacting dynamics,
stochastic averaging methods have proved useful for
vibro-impacting systems that experience a wide range of impact velocities \cite{DiMe79,RoSp86,FoBr96,DiIo04}.
If only low-velocity impacts are relevant, then it is useful to study (\ref{eq:N}).
In his PhD thesis \cite{Gr05}, Griffin studied (\ref{eq:N}) in the presence of additive white noise.
He found that noise blurs bifurcation diagrams and washes out high-period solutions
in the same manner as for smooth maps, such as the logistic map \cite{MaHa81,CrFa82}.
Recently it was shown that white noise added to (\ref{eq:f})
translates to additive white noise in (\ref{eq:N}), \cite{SiHo13}.
Such a noise formulation may be sensible for vibro-impacting systems
for which a forcing term or external fluctuations represent a significant source of uncertainty.

However, impact events themselves constitute a substantial source of randomness.
The purpose of this paper is to construct and analyse stochastic versions of (\ref{eq:N})
for which randomness stems purely from impact events.
We consider three different types of impact noise for (\ref{eq:f})
using an Ornstein-Uhlenbeck process with stationary density $\cN \left[ 0 ,\, \frac{\ee^2}{2 \nu} \right]$,
where $\ee \ll 1$ represents the noise amplitude and $\nu > 0$ is the correlation time.
We first consider uncertainties in $\Sigma$, then uncertainties in $f_R$, 
and lastly uncertainties in $f_R$ in the white noise (zero correlation time) limit.


Here let us indicate the forms of the stochastic maps that we obtain.
Coloured noise in $\Sigma$ leads to random perturbations in both the map and the switching condition as
\begin{equation}
\begin{bmatrix} x_{i+1} \\ y_{i+1} \end{bmatrix} = N_1(x_i,y_i) =
\begin{cases}
g(x_i,y_i) \;, & x_i + \kappa_1 \xi_i \le 0 \\
g(x_i,y_i) - \begin{bmatrix} \chi \sqrt{x_i + \kappa_1 \xi_i} \\ 0 \end{bmatrix} \;, & x_i + \kappa_1 \xi_i \ge 0
\end{cases} \;,
\label{eq:Na2}
\end{equation}
where $\xi_i \in \mathbb{R}$ are Gaussian random variables, and $\kappa_1 > 0$ is a constant.
For coloured noise in the impacting dynamics we obtain
\begin{equation}
\begin{bmatrix} x_{i+1} \\ y_{i+1} \end{bmatrix} = N_2(x_i,y_i) =
\begin{cases}
g(x_i,y_i) \;, & x_i \le 0 \\
g(x_i,y_i) - \begin{bmatrix} \chi \kappa_2(\xi_i) \sqrt{x_i} \\ 0 \end{bmatrix} \;, & x_i \ge 0
\end{cases} \;,
\label{eq:Nb2}
\end{equation}
for a particular nonlinear function $\kappa_2$.
For white noise ($\nu \to 0$) in the impacting dynamics the map takes the form
\begin{equation}
\begin{bmatrix} x_{i+1} \\ y_{i+1} \end{bmatrix} = N_3(x_i,y_i) =
\begin{cases}
g(x_i,y_i) \;, & x_i \le 0 \\
g(x_i,y_i) - \begin{bmatrix} \chi \kappa_3(r_i,h_i) \sqrt{x_i} \\ 0 \end{bmatrix}
+ \kappa_4(h_i) x_i \;, & x_i \ge 0
\end{cases} \;,
\label{eq:Nc2}
\end{equation}
where $r_i$ and $h_i$ are random variables,
and $\kappa_3 : \mathbb{R}^2 \to \mathbb{R}$ and $\kappa_4 : \mathbb{R} \to \mathbb{R}^2$
are nonlinear functions.
Notice that $N_1$ is stochastic for $x_i < 0$, whereas $N_2$ and $N_3$ are not.
This is because noise in $\Sigma$ generates anomalous crossings of $x_i = 0$.
Consequently $N_1$ exhibits stochastic dynamics for $\mu < 0$, while $N_2$ and $N_3$ do not.
$N_2$ and $N_3$ involve noise terms proportional to $\sqrt{x_i}$,
and for this reason exhibit increasing variability for larger values of $\mu > 0$.


To obtain a more detailed comparison of $N_1$, $N_2$ and $N_3$,
we first carefully derive (\ref{eq:N}) in \S\ref{sec:Nordmark}. 
In \S\ref{sec:oscillator} we introduce a prototypical compliant vibro-impacting system to illustrate our results.
In \S\ref{sec:randomness} we add randomness and derive (\ref{eq:Na2})-(\ref{eq:Nc2}).
Since (\ref{eq:Na2})-(\ref{eq:Nc2}) involve fundamentally different noise terms,
they exhibit different sensitivities to the noise amplitude $\ee$.
Therefore we use different values of $\ee$ in each of the different models,
in order to make appropriate comparisons.
For each $N_j$ ($j=1,2,3$) we write $\ee = \tilde{\ee}_j \alpha$
and identify the appropriate value $\tilde{\ee}_j$, so that
$N_1$, $N_2$ and $N_3$ display roughly the same dynamics for
the vibro-impacting system with the parameters of Fig.~\ref{fig:detBifDiag} and $\mu = 0.03$
(chosen for illustration) and $\alpha = 1$.
This enables us to quantitatively compare $N_1$, $N_2$ and $N_3$ in \S\ref{sec:dynamics}.
Conclusions are presented in \S\ref{sec:conc}.

\section{A derivation of the Nordmark map}
\label{sec:Nordmark}
\setcounter{equation}{0}

\begin{figure}[b!]
\begin{center}
\setlength{\unitlength}{1cm}
\begin{picture}(15,7.5)
\put(0,0){\includegraphics[height=7.5cm]{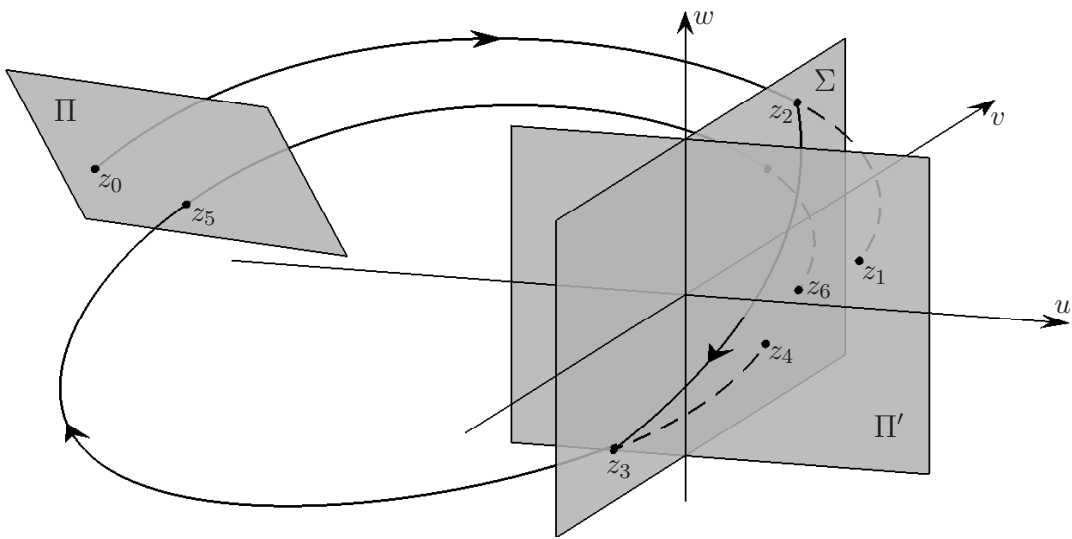}}
\put(14.5,3.25){\small $u$}
\put(13.65,5.75){\small $v$}
\put(9.7,7.1){\small $w$}
\put(1.2,5.8){$\Pi$}
\put(12.1,1.6){$\Pi'$}
\put(11.3,6.2){$\Sigma$}
\put(1.78,4.92){\small $z_0$}
\put(11.94,3.7){\small $z_1$}
\put(10.7,5.84){\small $z_2$}
\put(8.54,1.12){\small $z_3$}
\put(10.7,2.64){\small $z_4$}
\put(3.02,4.48){\small $z_5$}
\put(11.2,3.5){\small $z_6$}
\end{picture}
\caption{
A schematic diagram illustrating the construction of the Nordmark map (\ref{eq:N})
for the three-dimensional piecewise-smooth system (\ref{eq:f}).
The solid curve represents an orbit of (\ref{eq:f}).
The dashed curves show virtual extensions of this orbit into $u > 0$ as governed by $f_L$.
\label{fig:grazSchem3d}
}
\end{center}
\end{figure}

In this section we derive a return map for
the generic deterministic system (\ref{eq:f}) valid near the grazing bifurcation $\eta = 0$
and provide an explicit coordinate change that transforms the map to (\ref{eq:N}).
Such a derivation is given in \cite{DiBu08,DiBu01}.
We provide an explicit construction of (\ref{eq:N}) here
in order to provide a basis for deriving the stochastic maps in \S\ref{sec:randomness}.
Throughout this section we write $z = (u,v,w)$ for points of (\ref{eq:f}) in $\mathbb{R}^3$.

As discussed in \S\ref{sec:intro}, we assume that with $\eta = 0$ the system
(\ref{eq:f}) has a periodic orbit that intersects the origin, but is otherwise contained in $u < 0$.
This implies that $\dot{u} = 0$ for $f_L$ at the origin with $\eta = 0$
(i.e.~$e_1^{\sf T} f_L(0,0,0;0) = 0$).
For simplicity, we assume that we further have
\begin{equation}
e_1^{\sf T} f_L(u,0,w;\eta) = 0 \;,
\label{eq:PipCondition}
\end{equation}
for all $(u,w;\eta)$ in a neighbourhood of $(0,0;0)$,
which can usually be imposed by an appropriate coordinate change.
This assumption is particularly useful
in the case that randomness is present in the switching condition, see \S\ref{sub:s1}.

The key to deriving the Nordmark map (\ref{eq:N}) is selecting a convenient Poincar\'{e} section and
constructing a {\em discontinuity map} that accounts for the difference between
$f_L$ and $f_R$, that is, the difference between impacting and non-impacting dynamics.
We let $\Pi$ represent a generic Poincar\'{e} section of (\ref{eq:f})
that lies in $u < 0$ and intersects the grazing periodic orbit transversally,
and let $\Pi'$ be the coordinate plane $v=0$, see Fig.~\ref{fig:grazSchem3d}.

Given $z_0 \in \Pi$,
we define $z_1 \in \Pi'$ as the next intersection of the forward orbit
governed by $f_L$ (i.e.~ignoring the switching condition at $u=0$) with $\Pi'$.
If $u > 0$ for the point $z_1$, then $z_1$ does not represent the true intersection
of the orbit of (\ref{eq:f}) with $\Pi$.
Nevertheless, we study the return map on $\Pi'$
($z_1 \to z_6$) rather than the return map on $\Pi$ ($z_0 \to z_5$),
because the two maps are conjugate\removableFootnote{
Conjugacy between maps using two different Poincar\'{e} sections is discussed in, for instance, \cite{Me07}, pg.~156.
See also \cite{DiBu08}, pg.~280, for a discussion regarding $\Pi$ and $\Pi'$.
}
and the map on $\Pi'$ has a simpler form than the map on $\Pi$.

To derive the return map on $\Pi'$,
we must consider three additional points, $z_2$, $z_3$, and $z_4$, see Fig.~\ref{fig:grazSchem3d}.
The points $z_2$ and $z_3$ correspond to the entry and exit locations of the orbit with the impacting region, $u > 0$.
Then $z_1$ is obtained by travelling forward from the entry point $z_2$ to $\Pi'$ using $f_L$,
whilst $z_4$ is obtained by travelling backward from the exit point $z_3$ to $\Pi'$ using $f_L$.
The discontinuity map is defined as $z_4 = D(z_1)$.
If $u \le 0$ for the point $z_1$, then $z_4 = z_1$.

In order to derive an explicit expression for the discontinuity map, we consider the
three steps, $z_1 \to z_2$, $z_2 \to z_3$, and $z_3 \to z_4$, individually,
and expand $f_L$ and $f_R$ 
about $\eta = 0$ and the origin by writing
\begin{equation}
f_J(u,v,w;\eta) =
\begin{bmatrix}
\alpha_J v + \cO(2) \\
-\beta_J + \cO(1) \\
\gamma_J + \cO(1)
\end{bmatrix} \;,
\label{eq:alphabetagamma}
\end{equation}
where $\alpha_J, \beta_J, \gamma_J \in \mathbb{R}$, and $J = L,R$.
In (\ref{eq:alphabetagamma}), and throughout this section,
$\cO(k)$ is short-hand big-O notation for $\cO \left( \left( \sqrt{u},v,w,\eta \right)^k \right)$.
Notice that $\sqrt{u}$ is assumed to be of the same order as $v$, $w$, and $\eta$,
which is appropriate in view of the tangency between the grazing periodic orbit and $\Sigma$.
We assume
\begin{equation}
\alpha_L ,\, \alpha_R ,\, \beta_L ,\, \beta_R > 0 \;,
\label{eq:alphaLRbetaLRCondition}
\end{equation}
such that the tangency is quadratic for both $f_L$ and $f_R$,
and has the orientation depicted in Fig.~\ref{fig:grazSchem3d}.

From series expansions of the orbits governed by $f_L$ and $f_R$, see \cite{DiBu08,DiBu01},
we obtain the following formulas for the three steps in the discontinuity map in the case $u_1 > 0$
(writing $z_i = (u_i,v_i,w_i)$).
For $z_1 \to z_2$,
\begin{equation}
\begin{bmatrix} v_2 \\ w_2 \end{bmatrix} =
\begin{bmatrix}
\frac{\sqrt{2 \beta_L}}{\sqrt{\alpha_L}} \sqrt{u_1} + \cO(2) \\
w_1 - \frac{\sqrt{2} \gamma_L}{\sqrt{\alpha_L \beta_L}} \sqrt{u_1} + \cO(2)
\end{bmatrix} \;, \label{eq:v2w2}
\end{equation}
for $z_2 \to z_3$,
\begin{equation}
\begin{bmatrix} v_3 \\ w_3 \end{bmatrix} =
\begin{bmatrix}
-v_2 + \cO(2) \\
w_2 + \frac{2 \gamma_R}{\beta_R} v_2 + \cO(2)
\end{bmatrix} \;, \label{eq:v3w3}
\end{equation}
and for $z_3 \to z_4$,
\begin{equation}
\begin{bmatrix} u_4 \\ w_4 \end{bmatrix} =
\begin{bmatrix}
\frac{\alpha_L}{2 \beta_L} v_3^2 + \cO(3) \\
w_3 + \frac{\gamma_L}{\beta_L} v_3 + \cO(2)
\end{bmatrix} \;. \label{eq:u4w4}
\end{equation}
By combining (\ref{eq:v2w2})-(\ref{eq:u4w4})
we get $(u_4,w_4)$ in terms of $(u_1,w_1)$,
\begin{equation}
\begin{bmatrix} u_4 \\ w_4 \end{bmatrix} =
\begin{bmatrix}
u_1 + \cO(3) \\
w_1 - c \sqrt{u_1} + \cO(2)
\end{bmatrix} \;,
\label{eq:u4w42}
\end{equation}
where
\begin{equation}
c = \frac{2 \sqrt{2 \beta_L}}{\sqrt{\alpha_L}}
\left( \frac{\gamma_L}{\beta_L} - \frac{\gamma_R}{\beta_R} \right) \;.
\label{eq:c}
\end{equation}
The discontinuity map is then
\begin{equation}
D(u,w;\eta) =
\begin{cases}
\begin{bmatrix} u \\ w \end{bmatrix} \;, & u \le 0 \\
\begin{bmatrix} u + \cO(3) \\ w - c \sqrt{u} + \cO(2) \end{bmatrix} \;, & u \ge 0
\end{cases} \;.
\label{eq:D}
\end{equation}

To complete the map (\ref{eq:N}), which represents $z_1 \to z_6$ as shown in Fig.~\ref{fig:grazSchem3d},
we must combine $D(z_1)$ with the global return map $z_6 = G(z_4)$. 
$G$ depends on global properties of $f_L$ and is smooth, so we can write
\begin{equation}
G(u,w;\eta) = A \begin{bmatrix} u \\ w \end{bmatrix}
+ b \eta + \cO \left( (u,w,\eta)^2 \right) \;,
\label{eq:G}
\end{equation}
for some
\begin{equation}
A = \begin{bmatrix} a_{11} & a_{12} \\ a_{21} & a_{22} \end{bmatrix} \;, \qquad
b = \begin{bmatrix} b_1 \\ b_2 \end{bmatrix} \;,
\label{eq:Ab}
\end{equation}
where each $a_{ij}, b_i \in \mathbb{R}$.
Then the desired return map on $\Pi'$ is the composition $G \circ D$.

Finally we apply a coordinate change to convert the map to the normal form (\ref{eq:N})
that involves only three parameters, $\tau$, $\delta$ and $\chi$.
Under
\begin{equation}
\begin{bmatrix} x \\ y \\ \mu \end{bmatrix} =
\frac{1}{a_{12}^2 c^2} \begin{bmatrix}
1 & 0 & 0 \\
-a_{22} & a_{12} & b_1 \\
0 & 0 & (1-a_{22}) b_1 + a_{12} b_2
\end{bmatrix}
\begin{bmatrix} u \\ w \\ \eta \end{bmatrix} \;,
\label{eq:coordinateChange}
\end{equation}
and with higher order terms omitted,
$G \circ D$ transforms to (\ref{eq:N}) with
\begin{equation}
\tau = {\rm trace}(A) \;, \qquad
\delta = {\rm det}(A) \;, \qquad
\chi = {\rm sgn}(a_{12} c) \;.
\label{eq:taudeltachi}
\end{equation}
Since the Nordmark map keeps only leading order terms for $(u,w;\eta)$ near $(0,0;0)$,
terms that are linear in $x$ are omitted since they are dominated by $\sqrt{x}$,
refer to \cite{MoDe01} for a further discussion.
In (\ref{eq:coordinateChange}) we require $a_{12} \ne 0$ and $c \ne 0$.
These represent non-degeneracy conditions for the grazing bifurcation.

\section{An oscillator with compliant impacts}
\label{sec:oscillator}
\setcounter{equation}{0}

\begin{figure}[b!]
\begin{center}
\setlength{\unitlength}{1cm}
\begin{picture}(15,5)
\put(0,0){\includegraphics[height=5cm]{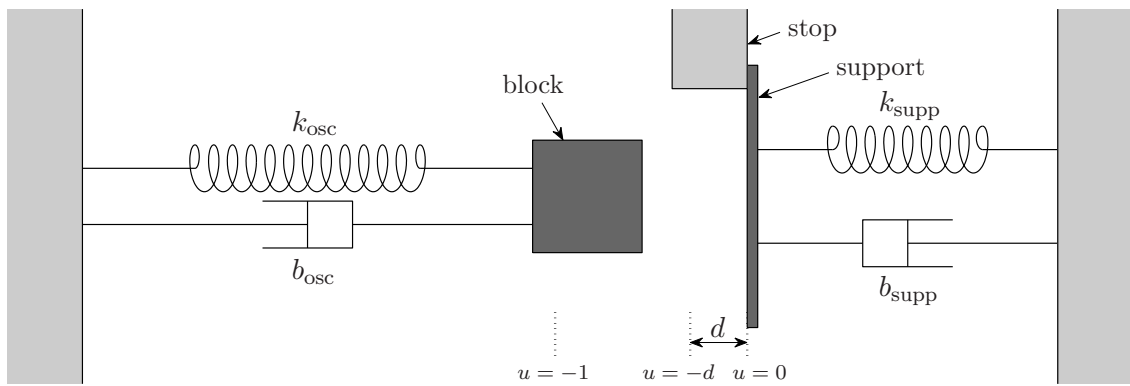}}
\put(6.78,.08){\scriptsize $u=-1$}
\put(8.44,.08){\scriptsize $u=-d$}
\put(9.64,.08){\scriptsize $u=0$}
\put(3.78,3.38){\small $k_{\rm osc}$}
\put(3.78,1.4){\small $b_{\rm osc}$}
\put(11.58,3.64){\small $k_{\rm supp}$}
\put(11.58,1.18){\small $b_{\rm supp}$}
\put(9.34,.63){\small $d$}
\put(6.58,3.88){\footnotesize block}
\put(11.02,4.12){\footnotesize support}
\put(10.38,4.64){\footnotesize stop}
\end{picture}
\caption{
A schematic diagram of the vibro-impacting system modelled by (\ref{eq:impactOscillator}).
This system exhibits a {\em regular} grazing bifurcation with a square-root singularity
because with $k_{\rm supp} > 0$ and $d > 0$
the equations of motion are discontinuous at the grazing point.
\label{fig:vibroImpactCompliant}
}
\end{center}
\end{figure}

To motivate and illustrate our results for stochastic Nordmark maps,
we consider the prototypical vibro-impacting system shown in Fig.~\ref{fig:vibroImpactCompliant}
and studied in \cite{SiHo13,MaAg06,MaIn08}.
This system consists of a harmonically forced one-degree-of-freedom linear oscillator that experiences
compliant (or soft) impacts with a support,
and we use the following non-dimensionalised equations to model the dynamics:
\begin{equation}
\ddot{u} = \begin{cases}
-k_{\rm osc} (u+1) - b_{\rm osc} \dot{u} + F \cos(t) \;, & u < 0 \\
-k_{\rm osc} (u+1) - (b_{\rm osc} + b_{\rm supp}) \dot{u} - k_{\rm supp} (u+d) + F \cos(t) \;, & u > 0
\end{cases} \;.
\label{eq:impactOscillator}
\end{equation}
Here $u(t)$ denotes the location of the block, which has the equilibrium position $u = -1$.
A rigid stop prevents the support reaching a position with $u < 0$,
and prestresses the support by a distance, $d > 0$.
The constants $k_{\rm osc}$, $b_{\rm osc}$, $k_{\rm supp}$ and $b_{\rm supp}$
represent non-dimensionalised spring and damping coefficients for the oscillator and support.
We neglect the mass of the support,
ignore energy loss at impacts,
and assume that whenever the block is not in contact with the support,
the support is located at $u=0$.
Experiments of simple vibro-impacting systems with compliant impacts
have shown that low-dimensional models such as (\ref{eq:impactOscillator})
can quantitatively match the physically observed dynamics near grazing bifurcations \cite{InPa06,InPa08b,InPa10}.

Here we treat the forcing amplitude, $F > 0$, as the primary bifurcation parameter,
and assume $k_{\rm osc}, b_{\rm osc} > 0$\removableFootnote{
If $k_{\rm osc} = 1$ the frequency of the steady-state and transient solutions are identical,
but this is not a problem because $b_{\rm osc} > 0$.
}.
The steady-state solution (behaviour in the limit $t \to \infty$) to (\ref{eq:impactOscillator}) with $u < 0$ is
\begin{equation}
u_{\rm ss}(t) = -1 + \frac{(k_{\rm osc}-1) \cos(t) + b_{\rm osc} \sin(t)}
{(k_{\rm osc}-1)^2 + b_{\rm osc}^2} F \;.
\end{equation}
When $F < F_{\rm graz}$, where
\begin{equation}
F_{\rm graz} = \sqrt{(k_{\rm osc}-1)^2 + b_{\rm osc}^2} \;,
\label{eq:Fgraz}
\end{equation}
the maximum value of $u_{\rm ss}(t)$ over one period is negative,
and so $u_{\rm ss}(t)$ is an attracting non-impacting periodic orbit of (\ref{eq:impactOscillator}).
The critical value, $F = F_{\rm graz}$, is a grazing bifurcation at which $u_{\rm ss}(t)$ has unit amplitude
and attains the value $u = 0$ at times $t = t_{\rm graz} + 2 \pi k$, for $k \in \mathbb{Z}$, where
\begin{equation}
t_{\rm graz} = {\rm tan}^{-1} \left( \frac{b_{\rm osc}}{k_{\rm osc} - 1} \right) \;,
\label{eq:tgraz}
\end{equation}
and $0 < t_{\rm graz} < \pi$\removableFootnote{
Since the coefficient of the $\sin(t)$ term of $u_{\rm ss}(t)$ is positive,
the maximum value of $u_{\rm ss}(t)$ must occur at some time with $0 < t_{\rm graz} < \pi$.
In addition, if $k_{\rm osc} > 1$, then $t_{\rm graz} < \frac{\pi}{2}$,
and if $0 < k_{\rm osc} < 1$, then $t_{\rm graz} > \frac{\pi}{2}$.
}.

To convert (\ref{eq:impactOscillator}) to the general form (\ref{eq:f}), we define
\begin{equation}
v = \dot{u} \;, \qquad
w = (t {\rm ~mod~} 2 \pi) - t_{\rm graz} \;, \qquad
\eta = F - F_{\rm graz} \;.
\label{eq:vw}
\end{equation}
Here the phase space of (\ref{eq:impactOscillator}) with (\ref{eq:vw}) is isomorphic to $\mathbb{R}^2 \times \mathbb{T}$,
rather than $\mathbb{R}^3$, but this does not affect the bifurcation structure near grazing.
For (\ref{eq:impactOscillator}) with (\ref{eq:vw}), the coefficients in (\ref{eq:alphabetagamma}),
which describe the behaviour of the system near the grazing point, are given by
\begin{equation}
\begin{gathered}
\alpha_L = 1 \;, \qquad
\beta_L = 1 \;, \qquad
\gamma_L = 1 \;, \\
\alpha_R = 1 \;, \qquad
\beta_R = 1 + k_{\rm supp} d \;, \qquad
\gamma_R = 1 \;,
\end{gathered}
\label{eq:ODECoeffsOsc}
\end{equation}
and by evaluating (\ref{eq:c}) with (\ref{eq:ODECoeffsOsc}) we obtain
\begin{equation}
c = \frac{2 \sqrt{2} k_{\rm supp} d}{1 + k_{\rm supp} d} \;.
\label{eq:cOsc}
\end{equation}
In addition, from the general solution to (\ref{eq:impactOscillator})
we find that the coefficients in the global map (\ref{eq:G}) are given by
\begin{equation}
A = {\rm exp} \left(
2 \pi \begin{bmatrix}
0 & 1 \\
-k_{\rm osc} & -b_{\rm osc}
\end{bmatrix} \right) \;, \qquad
b = \frac{1}{F_{\rm graz}}
\begin{bmatrix} 1 - a_{11} \\ -a_{21} \end{bmatrix} \;.
\end{equation}
\label{eq:mapCoeffOsc}

\section{Incorporating randomness into the Nordmark map}
\label{sec:randomness}
\setcounter{equation}{0}

To model noise and uncertainties we use the
one-dimensional Ornstein-Uhlenbeck process 
\begin{equation}
d\xi(t) = -\frac{1}{\nu} \xi(t) \,dt + \frac{\ee}{\nu} \,dW(t) \;,
\label{eq:xi}
\end{equation}
where $\ee, \nu > 0$ are constants
and $W(t)$ is a standard Brownian motion.
Given an initial value $\xi(0) = \xi_0$, 
at any positive time $\xi(t)$ is a Gaussian random variable with mean and variance
\begin{equation}
\mathbb{E} \left[ \xi(t) | \xi(0) = \xi_0 \right] =
\xi_0 \,{\rm e}^{\frac{-t}{\nu}} \;, \qquad
{\rm Var} \left[ \xi(t) | \xi(0) = \xi_0 \right] =
\frac{\ee^2}{2 \nu} \left( 1 - {\rm e}^{\frac{-2 t}{\nu}} \right) \;.
\label{eq:meanVarxi}
\end{equation}
In the limit $t \to \infty$,
$\xi(t) \sim \cN \left[ 0 ,\, \frac{\ee^2}{2 \nu} \right]$,
where $\cN[\mu,\sigma^2]$ denotes the Gaussian distribution of mean $\mu$ and variance $\sigma^2$.
The {\em correlation time} of (\ref{eq:xi}), defined as
$\int_0^\infty \frac{\mathbb{E} \left[ \xi(t) \xi(0) \right]}{{\rm Var} \left[ \xi(0) \right]} \,dt$,
with $\xi(0) \sim \cN \left[ 0 ,\, \frac{\ee^2}{2 \nu} \right]$,
is equal to $\nu$\removableFootnote{
Correlation time for coloured noise is defined at the start of \S 3 of \cite{HaJu95}.
This result follows from the formula
$\mathbb{E} \left[ \xi(t) \xi(s) \right] =
\frac{\ee^2}{2 \nu} \,{\rm e}^{\frac{-|t-s|}{\nu}}$.
}\removableFootnote{
With the alternate formulation,
\begin{equation}
d\zeta(t) = -\frac{1}{\nu} \zeta(t) \,dt + \frac{\sqrt{2} \ee}{\sqrt{\nu}} \,dW(t) \;,
\label{eq:zeta}
\end{equation}
the correlation time is again $\nu$
(since $\mathbb{E} \left[ \zeta(t) \zeta(s) \right] = \ee^2 \,{\rm e}^{\frac{-|t-s|}{\nu}}$), and
\begin{equation}
\mathbb{E} \left[ \zeta(t) | \zeta(0) = \zeta_0 \right] =
\zeta_0 \,{\rm e}^{\frac{-t}{\nu}} \;, \qquad
{\rm Var} \left[ \zeta(t) | \zeta(0) = \zeta_0 \right] =
\ee^2 \left( 1 - {\rm e}^{\frac{-2 t}{\nu}} \right) \;.
\end{equation}
In the limit $\nu \to 0$, the standard deviation of $\zeta(t)$
does not become unbounded, which makes (\ref{eq:zeta}) more useful than (\ref{eq:xi}) in some contexts.
}.

In our context, $\xi(t)$ is coloured noise
and the parameter $\ee$ governs the size of the noise.
Unlike white noise, $\xi(t)$ has an inherent time-scale, $\nu$, and is suitable
for modelling various types of uncertainties in mechanical systems,
such as background vibrations \cite{HaJu95}.
In the white noise limit,
forcing by $\xi(t)$ becomes a diffusion process $\ee \,dW(t)$\removableFootnote{
To prove this, we demonstrate that $\int_0^t \xi(s) \,ds \to \ee W(t)$.
By applying an integrating factor to (\ref{eq:xi}) we obtain
\begin{equation}
\xi(t) = \xi_0 {\rm e}^{\frac{-t}{\nu}} +
\frac{\ee}{\nu} \int_0^t {\rm e}^{\frac{u-t}{\nu}} \,dW(u) \;.
\end{equation}
By integrating this and reserving the order of the resulting double integral, we then obtain
\begin{equation}
\int_0^t \xi(s) \,ds = \xi_0 \nu \left( 1 - {\rm e}^{\frac{-t}{\nu}} \right) +
\ee \int_0^t 1 - {\rm e}^{\frac{u-t}{\nu}} \,dW(u) \;.
\end{equation}
Therefore $\int_0^t \xi(s) \,ds$ is Gaussian with mean,
$\xi_0 \nu \left( 1 - {\rm e}^{\frac{-t}{\nu}} \right) \to 0$ as $\nu \to 0$,
and variance, $\ee^2 \int_0^t \left( 1 - {\rm e}^{\frac{u-t}{\nu}} \right)^2 \,du \to \ee^2 t$ as $\nu \to 0$.
Hence, in the limit $\nu \to 0$, $\int_0^t \xi(s) \,ds$ has all properties of $\ee W(t)$
(mean $0$, variance $\ee^2 t$, correlation time, $0$).
}.


\subsection{Stochastic switching}
\label{sub:s1}

We first consider the following stochastic perturbation of (\ref{eq:f}),
\begin{equation}
\begin{bmatrix} \dot{u} \\ \dot{v} \\ \dot{w} \end{bmatrix} =
\begin{cases}
f_L(u,v,w;\eta) \;, & u + \xi(t) < 0 \\
f_R(u,v,w;\eta) \;, & u + \xi(t) > 0
\end{cases} \;,
\label{eq:fa}
\end{equation}
where $\xi(t)$ is given by (\ref{eq:xi}).
In (\ref{eq:fa}) randomness is present in the switching condition,
while evolution between switching events remains deterministic\removableFootnote{
If we want to take $\nu \to 0$ in (\ref{eq:fa}), we should set $\ee \propto \sqrt{\nu}$.
}.
We expect (\ref{eq:fa}) to be applicable to a wide variety of piecewise-smooth systems.
For the vibro-impacting system of \S\ref{sec:oscillator},
$\xi(t)$ may capture uncertainties in the point at which contact
between the block and support occurs or is lost.
For switched control systems, $\xi(t)$ may correspond to measurement errors
that produce variability in evaluations of switching rules \cite{Bo06b,Li03}.

Here we consider orbits of (\ref{eq:fa}) that are close to the grazing periodic orbit of (\ref{eq:f}).
Orbits of (\ref{eq:fa}) near grazing only spend short lengths of time in the region $u > 0$
while passing near the origin,
and for simplicity we suppose that the value of $\nu$ is large compared to such times.
In this case it is reasonable to approximate $\xi(t)$ by a constant while an orbit is near the origin.
With this approximation, the sum $u(t)+\xi(t)$ does not switch sign more than twice as the orbit passes near the origin,
which substantially simplifies our calculations below.

We let $\xi_i$ denote the value of $\xi(t)$ during the $i^{\rm th}$
instance that the orbit of (\ref{eq:fa}) passes near the origin.
The time between between consecutive traversals near the origin is well-approximated by the
period of the grazing periodic orbit, call it $T$.
With this approximation,
\begin{equation}
\xi_i \sim \cN \left[ \xi_{i-1} {\rm e}^{\frac{-T}{\nu}} ,\,
\frac{\ee^2}{2 \nu} \left( 1 - {\rm e}^{\frac{-2 T}{\nu}} \right) \right] \;.
\label{eq:xiiab}
\end{equation}

To derive the stochastic version of (\ref{eq:N}) for (\ref{eq:fa}) with (\ref{eq:xiiab}),
we first derive the induced stochastic discontinuity map.
Here condition (\ref{eq:PipCondition}) is useful,
as it implies that an orbit governed by $f_L$ attains a local maximum value of $u$ at an intersection with $\Pi'$. 
For $u_1 + \xi_i \le 0$, we conclude that $u(t) + \xi_i \le 0$ as the orbit passes near the origin, so that $u_4 = u_1$.
If instead $u_1 + \xi_i > 0$,
then the discontinuity map $D$ is given by (\ref{eq:u4w42})
except that $u + \xi_i$ appears inside the square root
because this quantity represents the distance from the switching condition.
That is,
\begin{equation}
\begin{bmatrix} u_4 \\ w_4 \end{bmatrix} =
\begin{bmatrix}
u_1 + \cO(3) \\
w_1 - c \sqrt{u_1 + \xi_i} + \cO(2)
\end{bmatrix} \;,
\label{eq:u4w42a}
\end{equation}
where $\cO(k) = \cO \left( \left( \sqrt{u},\sqrt{|\xi_i|},v,w,\eta \right)^k \right)$.
By combining (\ref{eq:u4w42a}) with the global map $G$,
applying the coordinate change (\ref{eq:coordinateChange}),
and dropping higher order terms, we obtain
\begin{equation}
\begin{bmatrix} x_{i+1} \\ y_{i+1} \end{bmatrix} = N_1(x_i,y_i) =
\begin{cases}
\begin{bmatrix} \tau & 1 \\ -\delta & 0 \end{bmatrix}
\begin{bmatrix} x_i \\ y_i \end{bmatrix} +
\begin{bmatrix} 0 \\ 1 \end{bmatrix} \mu \;, &
x_i + \frac{\xi_i}{a_{12}^2 c^2} \le 0 \\
\begin{bmatrix} \tau & 1 \\ -\delta & 0 \end{bmatrix}
\begin{bmatrix} x_i \\ y_i - \chi \sqrt{x_i + \frac{\xi_i}{a_{12}^2 c^2}} \end{bmatrix} +
\begin{bmatrix} 0 \\ 1 \end{bmatrix} \mu \;, &
x_i + \frac{\xi_i}{a_{12}^2 c^2} \ge 0
\end{cases} \;.
\label{eq:Na}
\end{equation}
$N_1$ is the stochastic Nordmark map corresponding to (\ref{eq:fa}).
Notice that randomness in the switching condition of (\ref{eq:fa})
has translated to randomness in both the switching condition of (\ref{eq:N})
and in the image of the map with $x_i > 0$.
In contrast, a piecewise-linear map
for which randomness is present purely in the switching condition is studied in \cite{Gl14b}.

In order to fairly compare $N_1$ with other stochastic versions of (\ref{eq:N}) in \S\ref{sec:dynamics},
we estimate the effective size of the stochastic contribution for our illustrative parameters values of the
prototypical system (\ref{eq:impactOscillator}) and a representative value of $\mu = 0.03$.
This motivates us to express $\ee$ in terms of a scaled parameter
$\ee = \tilde{\ee}_1 \alpha$, and to obtain comparable stochastic contributions for fixed $\alpha$ in the different cases.
The square-root term of $N_1$ is $\sqrt{x_i + \kappa_1 \xi_i}$,
where $\kappa_1 = \frac{1}{a_{12}^2 c^2}$.
If $x_i$ is large relative to $\kappa_1 \xi_i$,
then this term is well-approximated by $\sqrt{x_i} + \frac{\xi_i}{2 a_{12}^2 c^2 \sqrt{x_i}}$,
and so $\frac{\xi_i}{2 a_{12}^2 c^2 \sqrt{x_i}}$ estimates the additive stochastic contribution to $N_1$.
For the impact oscillator with the parameter values of Fig.~\ref{fig:detBifDiag}
(here $c = \sqrt{2}$ and $a_{12} \approx 0.1227$)
and using $x_i = 0.025$
(corresponding to the value of $x_i > 0$ in Fig.~\ref{fig:detBifDiag} for $\mu = 0.03$),
this quantity is approximately $100 \xi_i$.
With $\nu = 0.5$ (used in \S\ref{sec:dynamics}),
the standard derivation of the stochastic contribution is approximately $100 \ee$.
Therefore, for $\ee = \tilde{\ee}_1 \alpha$, where
\begin{equation}
\tilde{\ee}_1 = 0.0001 \;,
\label{eq:tildeee1}
\end{equation}
the standard deviation of the stochastic contribution is approximately $0.01$ when $\alpha = 1$.

\begin{figure}[b!]
\begin{center}
\setlength{\unitlength}{1cm}
\begin{picture}(6.8,5.1)
\put(0,0){\includegraphics[height=5.1cm]{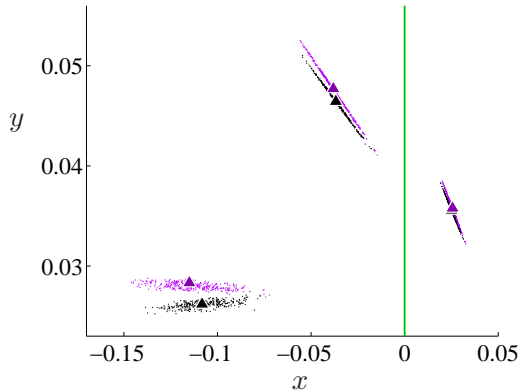}}
\put(3.76,0){\small $x$}
\put(0,3.5){\small $y$}
\end{picture}
\caption{
A phase portrait comparing the stochastic return map $N_1$, (\ref{eq:Na}), with a numerical solution to (\ref{eq:fa}).
The three groups of purple dots were obtained by numerically solving (\ref{eq:fa})
with (\ref{eq:xi}), $\nu = 0.5$ and $\ee = \tilde{\ee}_1$,
for the vibro-impacting system (\ref{eq:impactOscillator}) with (\ref{eq:vw}) and
$(k_{\rm osc},b_{\rm osc},k_{\rm supp},b_{\rm supp},d) = (4.5,0.3,10,0,0.1)$ (as in Fig.~\ref{fig:detBifDiag})
and $F \approx F_{\rm graz} + 0.005558$ (which corresponds to $\mu = 0.03$).
More precisely, $1000$ points on $\Pi'$ (labelled $z_1$ in Fig.~\ref{fig:grazSchem3d})
were obtained by numerically solving (\ref{eq:fa}),
and these were transformed to $(x,y)$-coordinates
by applying (\ref{eq:coordinateChange}) and (\ref{eq:vw})
to produce the purple dots.
The three groups of black dots are $1000$ iterates of $N_1$ with (\ref{eq:xiiab})
and parameter values chosen to match the vibro-impacting system
(specifically, $\nu = 0.5$, $\ee = \tilde{\ee}_1$, $T = 2 \pi$, $\mu = 0.03$,
$\tau \approx 0.5813$, $\delta \approx 0.1518$, $\chi = 1$, $c = \sqrt{2}$ and $a_{12} \approx 0.1227$).
The deterministic $3$-cycles of (\ref{eq:fa}) and (\ref{eq:Na}) are shown with triangles.
\label{fig:compareA_03}
}
\end{center}
\end{figure}

To illustrate the accuracy of $N_1$,
Fig.~\ref{fig:compareA_03} compares iterates of $N_1$ (black dots) with
a numerical solution to (\ref{eq:fa}) for the vibro-impacting system of \S\ref{sec:oscillator} (purple dots)
using $\ee = \tilde{\ee}_1$.
For the given parameter values, the system has an attracting $3$-cycle in the absence of noise.
For this reason, both sets of points are grouped about the $3$-cycle.
The two sets of points are slightly separated.
This is because the form of the deterministic Nordmark map does not include higher order terms of the true return map,
as observed by the separation of the values taken by the deterministic $3$-cycles shown in Fig.~\ref{fig:compareA_03}.
The size and shape of the spread of the two sets of randomly generated points are similar,
as is their location relative to the deterministic values of the map.
This demonstrates that $N_1$ can accurately capture the stochastic dynamics of (\ref{eq:fa}).
A more precise characterisation of the accuracy of $N_1$ is beyond the scope of this paper.

\subsection{Additive coloured noise with a large correlation time}
\label{sub:s2}

Next we consider the case where randomness and uncertainty in (\ref{eq:f}) is associated with $f_R$.
For mechanical systems with impacts,
this corresponds to variability in the evolution of the system during an impact. 
For simplicity we include noise in only the $v$-component of $f_R$, that is
\begin{equation}
\begin{bmatrix} \dot{u} \\ \dot{v} \\ \dot{w} \end{bmatrix} =
\begin{cases}
f_L(u,v,w;\eta) \;, & u < 0 \\
f_R(u,v,w;\eta) + \begin{bmatrix} 0 \\ \xi(t) \\ 0 \end{bmatrix} \;, & u > 0
\end{cases} \;,
\label{eq:fbc}
\end{equation}
where $\xi(t)$ is given by (\ref{eq:xi}).
Indeed, for the vibro-impacting system of Fig.~\ref{fig:vibroImpactCompliant},
if noise is incorporated into the force on the block when it is in contact with the support,
then the equations of motion may be put in the form (\ref{eq:fbc}).
With noise added to the $u$-component of $f_R$ (or $f_L$),
orbits may cross $u = 0$ many times in a short time frame
which makes the system substantially more difficult to analyse.
We leave such considerations for future work.

As in \S\ref{sub:s1}, we consider near-grazing orbits
and assume that the value of $\nu$ is much larger than the time each orbit spends in the region $u > 0$.
In this case it is reasonable to treat $\xi(t)$ as constant while $u > 0$.
During the $i^{\rm th}$ instance that an orbit passes near the origin,
we approximate $\xi(t)$ by $\xi_i$, distributed according to (\ref{eq:xiiab}).
In this scenario the three components of the discontinuity map (\ref{eq:v2w2})-(\ref{eq:u4w4})
are unchanged except that $\beta_R$ is replaced by $\beta_R - \xi_i$ in (\ref{eq:v3w3})
(because the $v$-component of the system with $u > 0$
is given by $-\beta_R + \xi_i + \cO(1)$, see (\ref{eq:alphabetagamma})).
By combining (\ref{eq:v2w2})-(\ref{eq:u4w4})
we find that for $u_1 > 0$ the discontinuity map is given by
\begin{equation}
\begin{bmatrix} u_4 \\ w_4 \end{bmatrix} =
\begin{bmatrix}
u_1 + \cO(3) \\
w_1 - c \left( \frac{\frac{\gamma_L}{\beta_L} - \frac{\gamma_R}{\beta_R - \xi_i}}
{\frac{\gamma_L}{\beta_L} - \frac{\gamma_R}{\beta_R}} \right) \sqrt{u_1} + \cO(2)
\end{bmatrix} \;,
\label{eq:u4w42b}
\end{equation}
and therefore the corresponding stochastic Nordmark map is
\begin{equation}
\begin{bmatrix} x_{i+1} \\ y_{i+1} \end{bmatrix} = N_2(x_i,y_i) =
\begin{cases}
\begin{bmatrix} \tau & 1 \\ -\delta & 0 \end{bmatrix}
\begin{bmatrix} x_i \\ y_i \end{bmatrix} +
\begin{bmatrix} 0 \\ 1 \end{bmatrix} \mu \;, & x_i \le 0 \\
\begin{bmatrix} \tau & 1 \\ -\delta & 0 \end{bmatrix}
\begin{bmatrix} x_i \\ y_i - \chi
\left( \frac{\frac{\gamma_L}{\beta_L} - \frac{\gamma_R}{\beta_R - \xi_i}}
{\frac{\gamma_L}{\beta_L} - \frac{\gamma_R}{\beta_R}} \right)
\sqrt{x_i} \end{bmatrix} +
\begin{bmatrix} 0 \\ 1 \end{bmatrix} \mu \;, & x_i \ge 0
\end{cases} \;.
\label{eq:Nb}
\end{equation}
Notice that with $\xi_i = 0$, $N_2$ is identical to (\ref{eq:N}).

We can write the stochastic component of $N_2$ as $\chi \kappa_2(\xi_i) \sqrt{x_i}$,
where $\kappa_2(\xi_i) = \frac{\frac{\gamma_L}{\beta_L} - \frac{\gamma_R}{\beta_R - \xi_i}}
{\frac{\gamma_L}{\beta_L} - \frac{\gamma_R}{\beta_R}}$.
With the parameter values of the impact oscillator (\ref{eq:ODECoeffsOsc}), and $k_{\rm supp} d = 1$, we have
$\kappa_2(\xi_i) = 2 - \frac{2}{2-\xi_i} \approx 1 - \frac{\xi_i}{2}$.
Therefore the noise provides a multiplicative stochastic contribution of approximately
$\frac{\sqrt{x_i} \xi_i}{2}$, ignoring signs.
In order to compare the effect of the noise to the other cases,
we write $\ee = \tilde{\ee}_2 \alpha$, choosing $\tilde{\ee}_2$
so that the standard deviation of the stochastic contribution is $0.01$ when $\alpha = 1$.
For $x_i = 0.025$ and $\nu = 0.5$ (as in \S\ref{sub:s1}),
the standard deviation of $\frac{\sqrt{x_i} \xi_i}{2}$ is approximately $0.08 \ee$, so we therefore choose
\begin{equation}
\tilde{\ee}_2 = 0.125 \;.
\label{eq:tildeee2}
\end{equation}

\begin{figure}[b!]
\begin{center}
\setlength{\unitlength}{1cm}
\begin{picture}(6.8,5.1)
\put(0,0){\includegraphics[height=5.1cm]{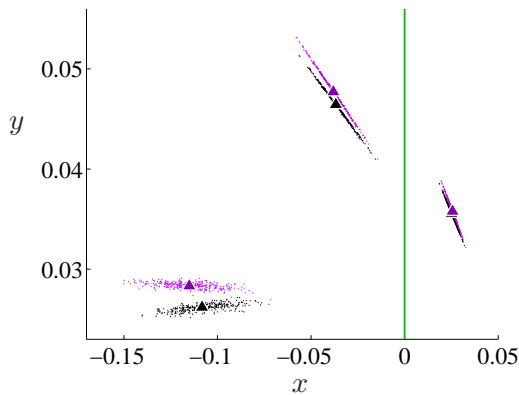}}
\put(3.76,0){\small $x$}
\put(0,3.5){\small $y$}
\end{picture}
\caption{
A phase portrait comparing the stochastic return map $N_2$ with a numerical solution to (\ref{eq:fbc}).
The three groups of purple dots were obtained by numerically solving (\ref{eq:fbc})
with (\ref{eq:xi}), $\nu = 0.5$ and $\ee = \tilde{\ee}_2$,
for the vibro-impacting system (\ref{eq:impactOscillator}) with (\ref{eq:vw})
using the same parameter values as in Fig.~\ref{fig:compareA_03}.
The three groups of black dots are $1000$ iterates of $N_2$ with (\ref{eq:xiiab})
and parameter values matching those of the vibro-impacting system
(refer to the caption of Fig.~\ref{fig:compareA_03}).
The deterministic $3$-cycles of (\ref{eq:fa}) and (\ref{eq:Na}) are shown with triangles.
\label{fig:compareB_03}
}
\end{center}
\end{figure}

Fig.~\ref{fig:compareB_03} compares iterates of $N_2$
to intersections with $\Pi'$ of a numerical solution to (\ref{eq:fbc}) using $\ee = \tilde{\ee}_2$.
As expected the two sets of points are
similarly distributed about the deterministic $3$-cycle.

\subsection{Additive coloured noise with a small correlation time}
\label{sub:s3}

Lastly we consider (\ref{eq:fbc}) in the white noise limit, $\nu = 0$.
In this case (\ref{eq:fbc}) reduces to a diffusion process forced by white noise,
specifically $\xi(t) \,dt$ is replaced by $\ee \,dW(t)$.
By using (\ref{eq:alphabetagamma}) to expand $f_R$,
(\ref{eq:fbc}) for $u > 0$ may be written as the three-dimensional stochastic differential equation
\begin{equation}
\begin{bmatrix} du(t) \\ dv(t) \\ dw(t) \end{bmatrix} =
\begin{bmatrix} \alpha_R v(t) + \cO(2) \\ -\beta_R + \cO(1) \\ \gamma_R + \cO(1) \end{bmatrix}
dt + \ee
\begin{bmatrix} 0 \\ 1 \\ 0 \end{bmatrix} dW(t) \;.
\label{eq:fRc}
\end{equation}	
To obtain a stochastic Nordmark map corresponding to this scenario,
we first derive the stochastic version of the middle component of the discontinuity map (\ref{eq:v3w3}).
Given an initial point $(u(0),v(0),w(0)) = (0,v_2,w_2)$, where $v_2 > 0$ and $w_2 \in \mathbb{R}$ are small,
the desired values of $v_3$ and $w_3$ are given by the point $(0,v_3,w_3)$ of {\em first return}
for the stochastic process (\ref{eq:fRc}) to $u = 0$.
First return or first passage problems are an important class of theoretical problems
in stochastic calculus with applications traditionally in finance and chemical kinetics \cite{GrVa99,Re01,Sc10}.

We approximate (\ref{eq:fRc}) by keeping only the leading order contributions, i.e.
\begin{equation}
\begin{bmatrix} du(t) \\ dv(t) \end{bmatrix} =
\begin{bmatrix} \alpha_R v(t) \\ -\beta_R \end{bmatrix}
dt + \ee
\begin{bmatrix} 0 \\ 1 \end{bmatrix} dW(t) \;, \qquad
\begin{bmatrix} u(0) \\ v(0) \end{bmatrix} =
\begin{bmatrix} 0 \\ v_2 \end{bmatrix} \;,
\label{eq:fRc2}
\end{equation}
together with $w(t) = w_2 + \gamma_R t$.
With this approximation we are able to provide an explicit expression
for the joint probability density function of the return location and time. 
A formal justification for the omission of the higher order terms is left for future work.

We introduce the change of variables\removableFootnote{
This choice of scaling leads to the simplest expression for the noisy Nordmark map
because here $w_3$ has an explicit linear dependence on $v_2$.
If we scale such that the parameter is in the drift, as in \cite{AtCl94},
then $r \propto \frac{1}{v_2}$ and we get $w_3 = w_2 + \frac{\gamma_R v_2^2 r}{\ee^2}$,
i.e.~$w_3 - w_2$ appears proportional to $v_2^2$.
If we instead scale such that the parameter is in $q(0)$,
then $r \propto v_2$ and we get $w_3 = w_2 + \frac{\gamma_R \ee^2 r}{\beta_R^2}$,
i.e.~$w_3 - w_2$ appears independent of $v_2$.
}
\begin{equation}
p = \frac{\beta_R}{\alpha_R v_2^2} u \;, \qquad
q = \frac{1}{v_2} v \;, \qquad
s = \frac{\beta_R}{v_2} t \;.
\label{eq:pqs}
\end{equation}
This puts (\ref{eq:fRc2}) in a standard form studied in \cite{AtCl94}\removableFootnote{
It is tempting to write $\varrho$ instead of $\sqrt{\varrho}$,
but this leads to powers of $\frac{1}{4}$ in (\ref{eq:varrho2}), which is not so nice.
},
\begin{equation}
\begin{bmatrix} dp(s) \\ dq(s) \end{bmatrix} =
\begin{bmatrix} q(s) \\ -1 \end{bmatrix} ds +
\sqrt{\varrho} \begin{bmatrix} 0 \\ 1 \end{bmatrix} dW(s) \;, \qquad
\begin{bmatrix} p(0) \\ q(0) \end{bmatrix} =
\begin{bmatrix} 0 \\ 1 \end{bmatrix} \;,
\label{eq:dpdq}
\end{equation}
where
\begin{equation}
\varrho = \frac{\ee^2}{\beta_R v_2} \;.
\label{eq:varrho}
\end{equation}
Specifically, $q$ is a diffusion process with constant drift,
and $p(s) = \int_0^s q(\tilde{s}) \,d\tilde{s}$.
Therefore $p(s)$ may be interpreted as {\em integrated Brownian motion with constant drift}.

We let $\cF(r,h;\varrho)$ denote the joint probability density function
for the first return of (\ref{eq:dpdq}) to $p = 0$, at a time $s = r > 0$, and location $q = -h < 0$.
In \cite{Mc63}, McKean derived an explicit expression for $\cF$ in the case of zero drift
by computing the inverse Kontorovich-Lebedev transform\removableFootnote{
This transformation uses a modified Bessel function of the second kind
as the kernel \cite{Er54,Ya96}.
}
of the corresponding renewal equation.
In \cite{AtCl94}, Atkinson and Clifford extended this result to the case of non-zero drift
by applying the Radon-Nikodyn derivative.
Specifically
\begin{equation}
\cF(r,h;\varrho) = \frac{\sqrt{3} h}{\pi \varrho r^2}
\,{\rm exp} \left( \frac{-1}{2 \varrho r}
\left[ (r-2)^2 - 2 (r-2)(h-1) + 4 (h-1)^2 \right] \right)
{\rm erf} \left( \frac{\sqrt{6 h}}{\sqrt{\varrho r}} \right) \;,
\label{eq:F}
\end{equation}
where ${\rm erf}(\cdot)$ is the error function.
The constant $\varrho > 0$ governs the shape of $\cF$.
The limit $\varrho \to 0$ corresponds to the deterministic case, for which  $r = 2$ and $h = 1$.
With a small value of $\varrho$, $\cF$ is roughly Gaussian.
In contrast, the limit $\varrho \to \infty$ corresponds to the case of no drift, as in \cite{Mc63},
or to the limit $v_2 \to 0$. 
In this limit the marginal probability density function for $r$
is asymptotically proportional to $r^{-\frac{5}{4}}$, for large $r$,
and so is long-tailed \cite{AtCl94}.

In view of the scaling (\ref{eq:pqs}),
the stochastic version of (\ref{eq:v3w3}) corresponding to (\ref{eq:fRc}) is given by\removableFootnote{
We treat $r$ and $h$ as $\cO(1)$ quantities.
}
\begin{equation}
\begin{bmatrix} v_3 \\ w_3 \end{bmatrix} =
\begin{bmatrix}
-h_i v_2 + \cO(2) \\
w_2 + \frac{\gamma_R r_i}{\beta_R} v_2 + \cO(2)
\end{bmatrix} \;,
\label{eq:v3w3c}
\end{equation}
where $r_i$ and $h_i$ have the joint probability density function (\ref{eq:F}).
By combining (\ref{eq:v3w3c}) with (\ref{eq:v2w2}) and (\ref{eq:u4w4}) we obtain
\begin{equation}
\begin{bmatrix} u_4 \\ w_4 \end{bmatrix} =
\begin{bmatrix}
h_i^2 u_1 + \cO(3) \\
w_1 - c \left( \frac{\frac{\gamma_L(h_i+1)}{2 \beta_L}-\frac{\gamma_R r_i}{2 \beta_R}}
{\frac{\gamma_L}{\beta_L}-\frac{\gamma_R}{\beta_R}} \right) \sqrt{u_1} + \cO(2)
\end{bmatrix} \;,
\label{eq:u4w4c}
\end{equation}
which represents the stochastic version of the discontinuity map for points with $u_1 > 0$.
Then using (\ref{eq:u4w4c}) we arrive at the following stochastic Nordmark map\removableFootnote{
I may wish to rewrite some of my Matlab code to match the notation used here
(e.g.~$\alpha_L$ instead of $\beta_L$, and $\varrho$ instead of $-\frac{1}{\kappa}$).
}
\begin{equation}
\begin{bmatrix} x_{i+1} \\ y_{i+1} \end{bmatrix} = N_3(x_i,y_i) =
\begin{cases}
\begin{bmatrix} \tau & 1 \\ -\delta & 0 \end{bmatrix}
\begin{bmatrix} x_i \\ y_i \end{bmatrix} +
\begin{bmatrix} 0 \\ 1 \end{bmatrix} \mu \;, & x_i \le 0 \\
\begin{bmatrix}
\tau + a_{11} (h_i^2-1) & 1 \\
-\delta h_i^2 & 0 \end{bmatrix}
\begin{bmatrix} x_i \\
y_i - \chi \left( \frac{\frac{\gamma_L(h_i+1)}{2 \beta_L}-\frac{\gamma_R r_i}{2 \beta_R}}
{\frac{\gamma_L}{\beta_L}-\frac{\gamma_R}{\beta_R}} \right) \sqrt{x_i}
\end{bmatrix} +
\begin{bmatrix} 0 \\ 1 \end{bmatrix} \mu \;, & x_i \ge 0
\end{cases} \;.
\label{eq:Nc}
\end{equation}
In (\ref{eq:Nc}) we treat each pair $(r_i,h_i)$ as a
two-dimensional stochastic random variable with probability density function
$\cF(r_i,h_i;\varrho_i)$, where
\begin{equation}
\varrho_i = \frac{\ee^2 \sqrt{\alpha_L}}{\beta_R \sqrt{2 \beta_L} |a_{12} c| \sqrt{x}} \;,
\label{eq:varrho2}
\end{equation}
which results from combining (\ref{eq:varrho}) with
$v_2 \approx \frac{\sqrt{2 \beta_L}}{\sqrt{\alpha_L}} \sqrt{u_1}$, (\ref{eq:v2w2}),
and $u_1 = a_{12}^2 c^2 x$, (\ref{eq:coordinateChange}).

We now estimate the size of the stochastic contribution in $N_3$.
The leading order stochastic component of $N_3$ is $\chi \kappa_3(r_i,h_i) \sqrt{x_i}$,
where $\kappa_3(r_i,h_i) = \frac{\frac{\gamma_L(h_i+1)}{2 \beta_L}-\frac{\gamma_R r_i}{2 \beta_R}}
{\frac{\gamma_L}{\beta_L}-\frac{\gamma_R}{\beta_R}}$.
With (\ref{eq:ODECoeffsOsc}) and $k_{\rm supp} d = 1$, we can write
$\kappa_3(r_i,h_i) \sqrt{x_i} = \sqrt{x_i} + \left( h_i-1 - \frac{r_i-2}{2} \right) \sqrt{x_i}$,
where the second term represents the multiplicative stochastic contribution, ignoring the sign of $\chi$,
because the deterministic values of $h_i$ and $r_i$ are $1$ and $2$ respectively.

With a small value of $\varrho$ ($\varrho < 0.03$ is suitable),
$\cF(r,h;\varrho)$ is approximately Gaussian
because the effective noise amplitude in (\ref{eq:dpdq}) is small.
By (\ref{eq:varrho2}), this approximation is valid when, roughly speaking,
$\ee$ is not too large and $x_i$ is not too small.
From (\ref{eq:F}) we determine the covariance matrix of the Gaussian approximation to be
\begin{equation}
{\rm Cov}(r,h;\varrho) =
\frac{2 \varrho}{3} \begin{bmatrix} 4 & 1 \\ 1 & 1 \end{bmatrix} \;,
\label{eq:Theta}
\end{equation}
and it follows that in this approximation the linear combination
$h - \frac{r}{2}$ has standard deviation $\sqrt{\frac{2 \varrho}{3}}$.
The standard deviation of the stochastic contribution in $N_3$ is therefore approximately
$\sqrt{\frac{2 \varrho_i}{3}} \sqrt{x_i}$.
Following the previous cases, we write $\ee = \tilde{\ee}_3 \alpha$,
and choose $\tilde{\ee}_3$ so that the standard deviation of the stochastic contribution is $0.01$ when $\alpha = 1$. 
For $N_3$, this quantity is $\sqrt{\frac{2 \varrho_i}{3}} \sqrt{x_i} = x_i^{\frac{1}{4}} \ee \approx 0.46 \ee$,
using (\ref{eq:varrho2}), the parameter values from Fig.~\ref{fig:detBifDiag}, and $x_i = 0.025$.
Therefore we let
\begin{equation}
\tilde{\ee}_3 = 0.022 \;.
\label{eq:tildeee3}
\end{equation}
Here $\varrho_i \approx 0.008$ when $\alpha = 1$,
and so for these values the Gaussian approximation to (\ref{eq:F}) is justified.

\begin{figure}[t!]
\begin{center}
\setlength{\unitlength}{1cm}
\begin{picture}(6.8,5.1)
\put(0,0){\includegraphics[height=5.1cm]{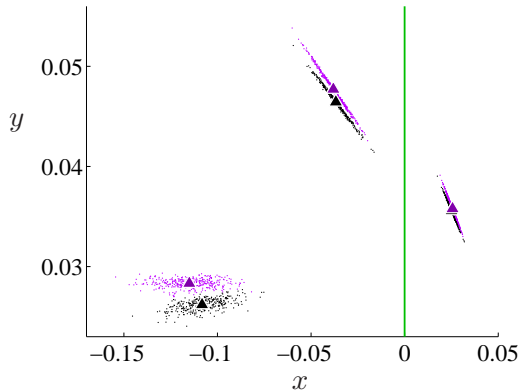}}
\put(3.76,0){\small $x$}
\put(0,3.5){\small $y$}
\end{picture}
\caption{
A phase portrait comparing the stochastic return map $N_3$ with a numerical solution to (\ref{eq:fbc}).
The three groups of purple dots were obtained by numerically solving (\ref{eq:fbc})
with (\ref{eq:xi}), $\nu = 0$
(in which case $\xi(t) \,dt$ is replaced by $\ee\,dW(t)$) and $\ee = \tilde{\ee}_3$,
for the vibro-impacting system (\ref{eq:impactOscillator}) with (\ref{eq:vw})
using the same parameter values as in Fig.~\ref{fig:compareA_03}.
The three groups of black dots are $1000$ iterates of $N_3$ with (\ref{eq:F}) and (\ref{eq:varrho2})
and parameter values matching those of the vibro-impacting system.
The deterministic $3$-cycles of (\ref{eq:fa}) and (\ref{eq:Na}) are shown with triangles.
\label{fig:compareC_03}
}
\end{center}
\end{figure}

Fig.~\ref{fig:compareC_03} compares iterates of $N_3$
to intersections with $\Pi'$ of a numerical solution to (\ref{eq:fbc}) using $\ee = \tilde{\ee}_3$.
As with the previous two figures,
this shows that $N_3$ can accurately capture the stochastic dynamics of (\ref{eq:fbc}).

\section{Stochastic dynamics}
\label{sec:dynamics}
\setcounter{equation}{0}

In this section we explore the dynamics of the three stochastic Nordmark maps,
$N_1$, $N_2$ and $N_3$, 
and discuss how the different forms of these maps is evident in their dynamical behaviour.
To briefly summarise,
$N_1$ applies to the system with stochastic switching (\ref{eq:fa}),
whereas $N_2$ and $N_3$ apply to (\ref{eq:fbc}).
For $N_1$ and $N_2$ it is assumed that the value of $\nu$
(the correlation time of the noise) is large relative to the times that orbits spend in $u > 0$.
As these times are rarely larger than
$t = 0.05$ for the parameter values considered here,
we take $\nu = 0.5$ in $N_1$ and $N_2$ to ensure that the correlation time is large enough.
The values of $\xi_i$ in $N_1$ and $N_2$ are distributed according to (\ref{eq:xiiab}).
$N_3$ corresponds to the limit $\nu \to 0$.
In $N_3$, $r_i$ and $h_i$ are distributed according to (\ref{eq:F})
and depend on the value of $\varrho_i$ (\ref{eq:varrho2}).

For each $N_j$ we have written
\begin{equation}
\ee = \tilde{\ee}_j \alpha \;,
\label{eq:tildeeej}
\end{equation}
where the $\tilde{\ee}_j$ are given by
(\ref{eq:tildeee1}), (\ref{eq:tildeee2}) and (\ref{eq:tildeee3}).
These values have been chosen such that
for a given value of $\alpha$, the size of the 
stochastic contribution in $N_1$, $N_2$ and $N_3$
is roughly the same, at least when $\mu = 0.03$.
For $\mu$ close to $0$, the stochastic contributions are noticeably different for these
choices of $\tilde{\ee}_j$.


In \S\ref{sub:muDependence} we look at stochastic bifurcation diagrams
in order to obtain a basic understanding of how the stochastic dynamics differs with the value of $\mu$.
In the subsequent parts of this section we study two-dimensional invariant densities
in order to gain a deeper understanding of the dynamics.

\subsection{The dependence of $\mu$ on the size of noise response}
\label{sub:muDependence}

\begin{figure}[b!]
\begin{center}
\setlength{\unitlength}{1cm}
\begin{picture}(12,18.4)
\put(0,12.4){\includegraphics[height=6cm]{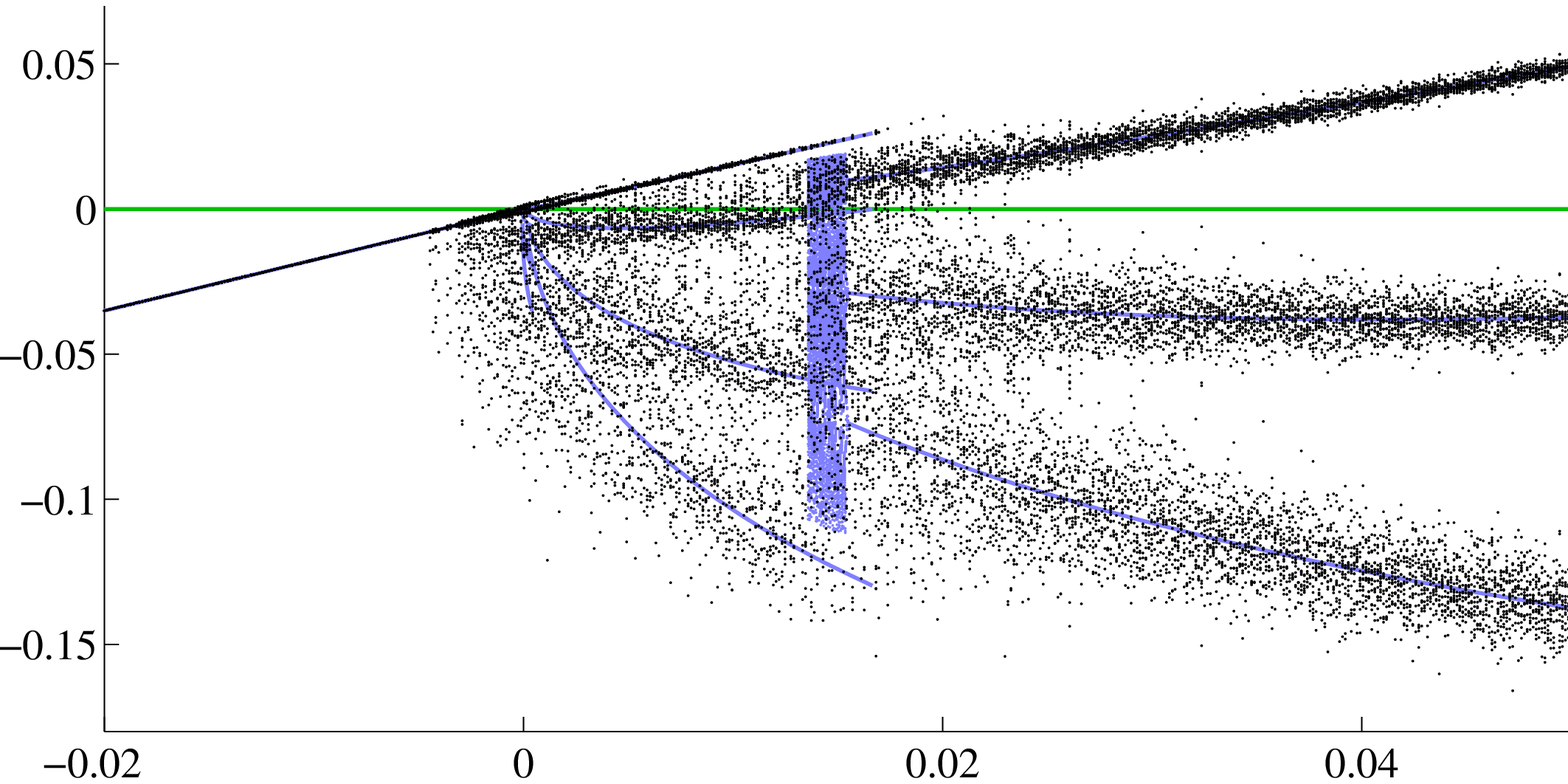}}
\put(0,6.2){\includegraphics[height=6cm]{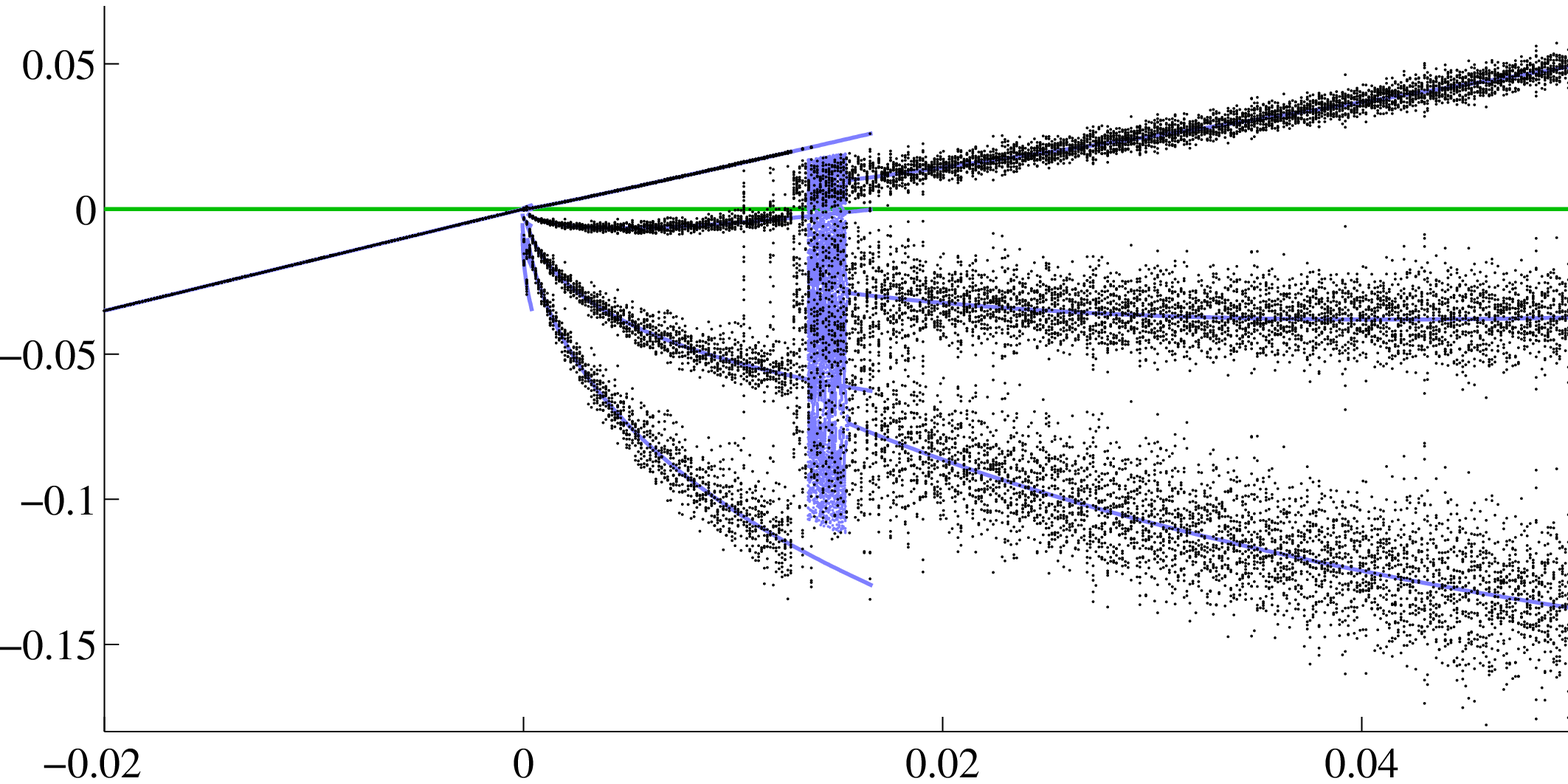}}
\put(0,0){\includegraphics[height=6cm]{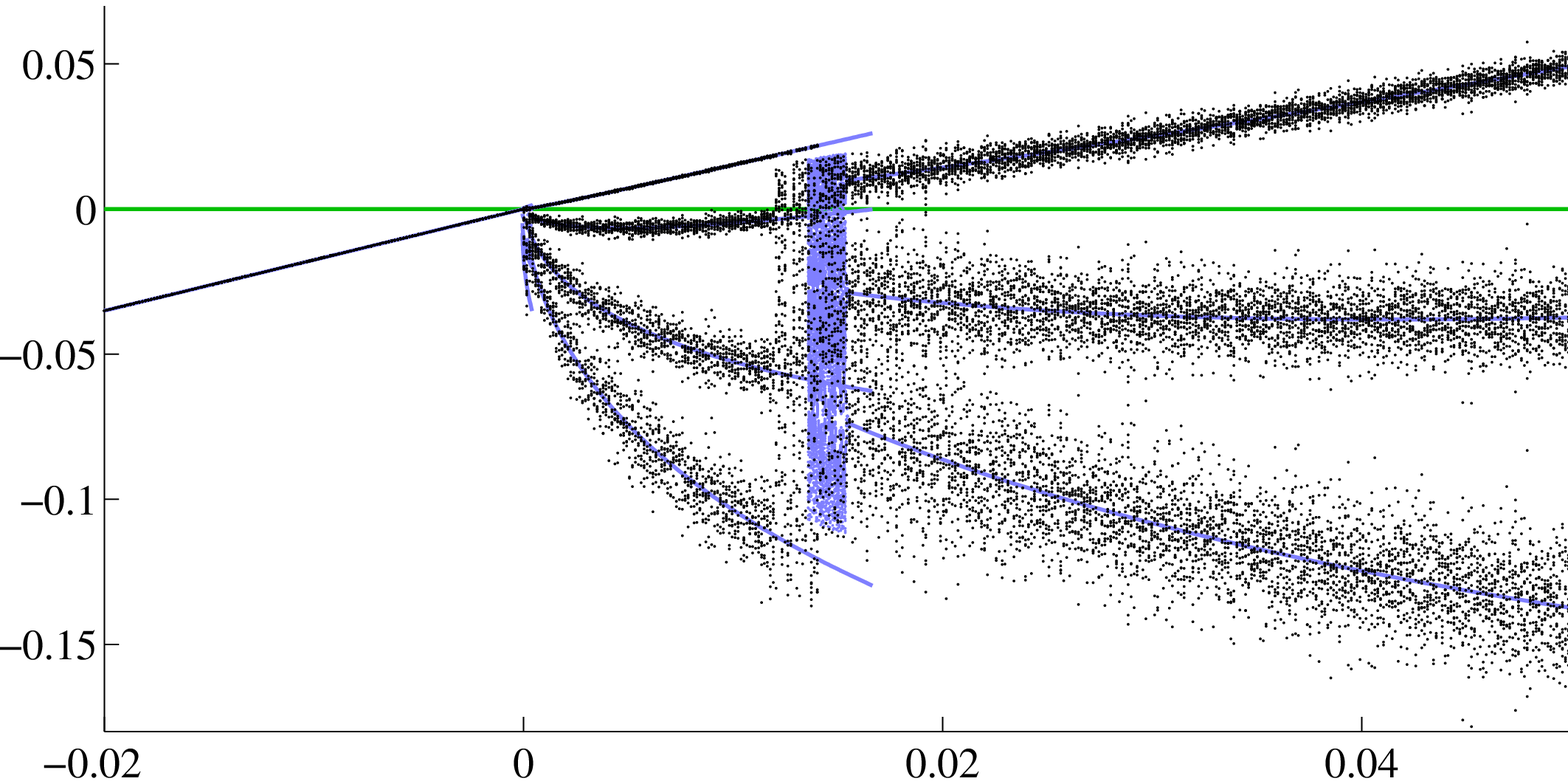}}
\put(6.2,12.4){\small $\mu$}
\put(6.2,6.2){\small $\mu$}
\put(6.2,0){\small $\mu$}
\put(0,16.2){\small $x$}
\put(0,10){\small $x$}
\put(0,3.8){\small $x$}
\put(1.7,18.1){\large \sf \bfseries A}
\put(1.7,11.9){\large \sf \bfseries B}
\put(1.7,5.7){\large \sf \bfseries C}
\end{picture}
\caption{
Bifurcation diagrams of the stochastic Nordmark maps
$N_1$ (panel A),
$N_2$ (panel B),
and $N_3$ (panel C),
with $\tau \approx 0.5813$, $\delta \approx 0.1518$, and $\chi = 1$.
These parameter values correspond to the vibro-impacting system (\ref{eq:impactOscillator}),
with $(k_{\rm osc},b_{\rm osc},k_{\rm supp},b_{\rm supp},d) = (4.5,0.3,10,0,0.1)$.
The black dots are iterates of $N_1$, $N_2$ and $N_3$ with transient points omitted.
The noise amplitudes are given by (\ref{eq:tildeee1}), (\ref{eq:tildeee2}) and (\ref{eq:tildeee3}) with $\alpha = 1$,
and $\nu = 0.5$ for panels A and B (panel C corresponds to $\nu = 0$).
In each panel the deterministic bifurcation diagram (Fig.~\ref{fig:detBifDiag}) is shown in blue.
\label{fig:noisyBifDiag}
}
\end{center}
\end{figure}

Fig.~\ref{fig:noisyBifDiag} shows stochastic versions of the bifurcation diagram shown in Fig.~\ref{fig:detBifDiag}
for the three stochastic Nordmark maps.
As expected, the noise blurs the bifurcation diagram.
In panel A, which corresponds to the map $N_1$,
for values of $\mu$ very close to zero (say $|\mu| < 0.002$)
the points are relatively highly spread.
This is because here the deterministic map has an attractor near $x=0$,
so in $N_1$ the sign of $x_i$ is often different to the sign of $x_i + \frac{\xi_i}{a_{12}^2 c^2}$.
That is, the choice of the half-map of $N_1$ is regularly determined by $\xi_i$ rather than $x_i$.
Furthermore, as shown in \S\ref{sub:s1}, the leading order component of the
stochastic contribution to the right half-map of $N_1$
is inversely proportional to $\sqrt{x_i}$,
which for very small $\mu$ is large relative to its value for $\mu$ away from zero.
In contrast, with $0.03 \le \mu \le 0.05$ say,
the underlying attracting $3$-cycle is sufficiently far from $x=0$
so that the sign of $x_i$ rarely differs from that of $x_i + \frac{\xi_i}{a_{12}^2 c^2}$.
The points are randomly distributed\hfill

\clearpage																						

\noindent
about the $3$-cycle
due to noise in the right half-map of $N_1$.
The size of the deviation decreases with increasing $\mu$,
because the strength of attraction of the $3$-cycle increases with $\mu$.

In panels B and C, which correspond to $N_2$ and $N_3$ respectively,
the bifurcation diagrams show no variability for $\mu < 0$.
This is because for $\mu < 0$, $(x^L,y^L)$ (\ref{eq:xLyL}) is a fixed point of $N_2$ and $N_3$.
For $0.002 \le \mu \le 0.01$,
the size of deviations about the underlying attracting $4$-cycle increases with $\mu$.
This is primarily because the coefficient of the $\sqrt{x}$-term of $N_2$ and $N_3$ is random,
and for the $4$-cycle this value of $x$ increases with $\mu$.
For $0.03 \le \mu \le 0.05$,
the size of deviations varies little with $\mu$ because
the increased variability caused by a larger value of $x$ 
is balanced by the fact that the strength of attraction of the $3$-cycle increases.
Panels B and C are similar,
suggesting that the value of $\nu$ has little effect on the long-term dynamics,
although panel C shows slightly more variability for very small values of $\mu > 0$.

\subsection{Invariant densities about an attracting periodic solution}

\begin{figure}[b!]
\begin{center}
\setlength{\unitlength}{1cm}
\begin{picture}(13.8,10.2)
\put(0,4.7){\includegraphics[height=5.1cm]{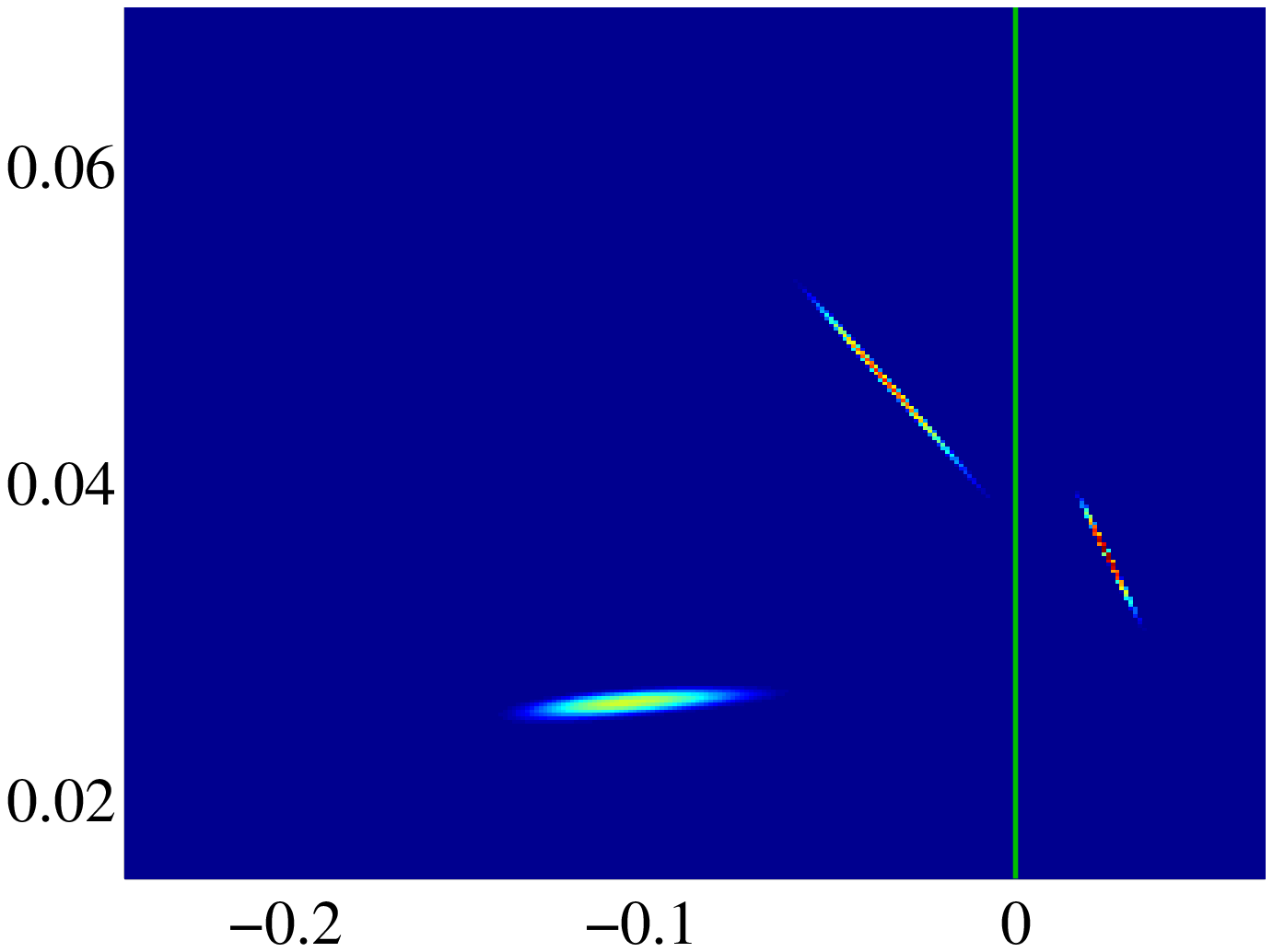}}
\put(7,4.7){\includegraphics[height=5.1cm]{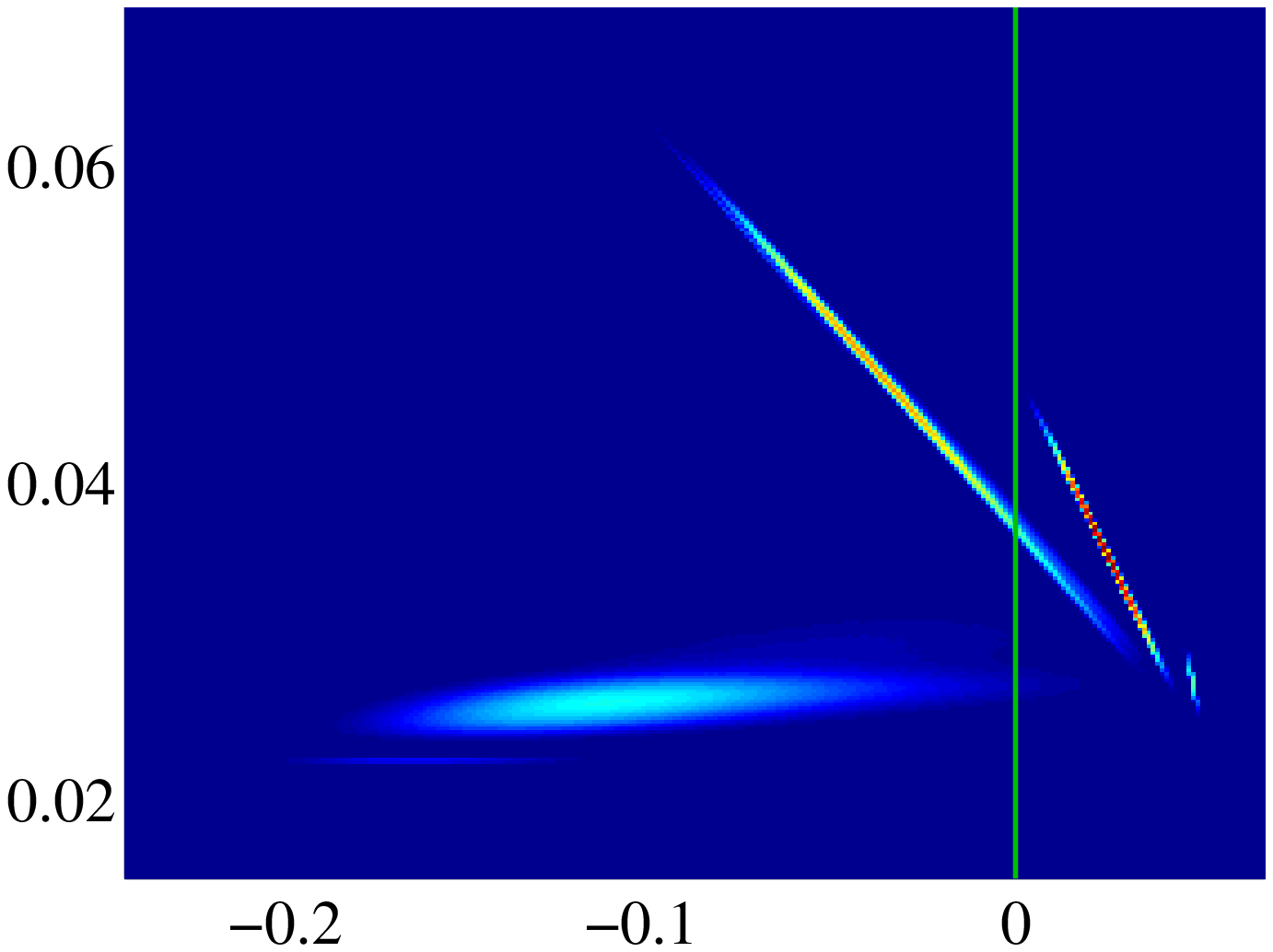}}
\put(2.4,0){\includegraphics[height=4.5cm]{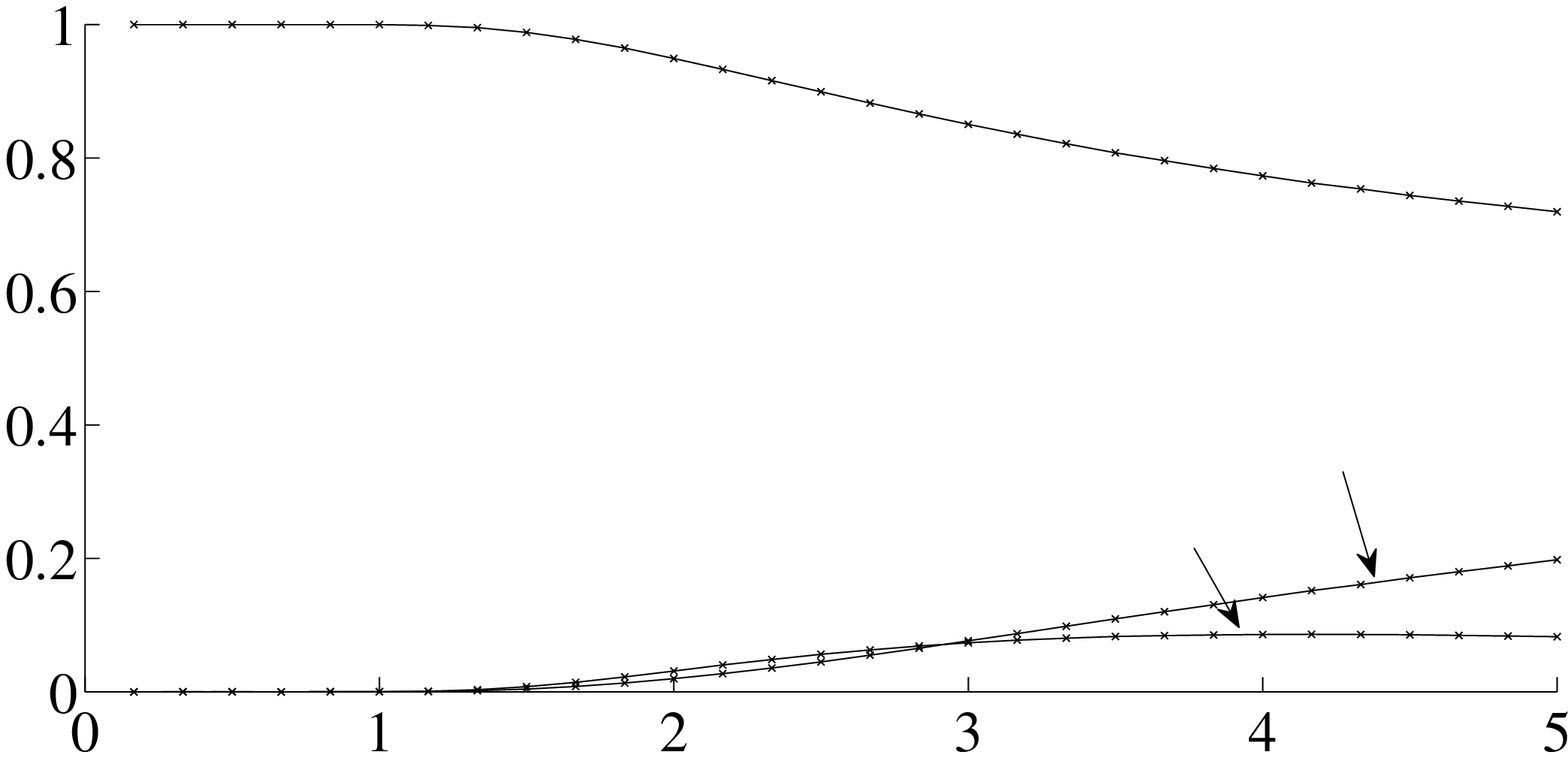}}
\put(3.42,10){\small $\alpha = 1$}
\put(10.42,10){\small $\alpha = 3$}
\put(3.76,4.7){\small $x$}
\put(.3,8.2){\small $y$}
\put(10.76,4.7){\small $x$}
\put(7.3,8.2){\small $y$}
\put(7.2,0){\small $\alpha$}
\put(8,3.44){\small $\sigma_3$}
\put(8.87,1.45){\small $\sigma_2$}
\put(9.33,2){\small $1 - \sigma_2 - \sigma_3$}
\put(.5,9.9){\large \sf \bfseries A}
\put(7.5,9.9){\large \sf \bfseries B}
\put(2.4,4.2){\large \sf \bfseries C}
\end{picture}
\caption{
Panel A shows the invariant density of $N_1$
with the same parameter values as Fig.~\ref{fig:noisyBifDiag} and $\mu = 0.03$.
The value of the density is indicated by colour
(dark red -- the maximum value of the density; dark blue -- zero).
Panel B shows the invariant density at three times the noise amplitude as panel A.
Panel C plots $\sigma_j$ -- the fraction of instances that points return to $x>0$
in $j$ iterations (\ref{eq:sigma}) -- against the noise amplitude.
\label{fig:invDensityA_03}
}
\end{center}
\end{figure}

\begin{figure}[b!]
\begin{center}
\setlength{\unitlength}{1cm}
\begin{picture}(13.8,10.2)
\put(0,4.7){\includegraphics[height=5.1cm]{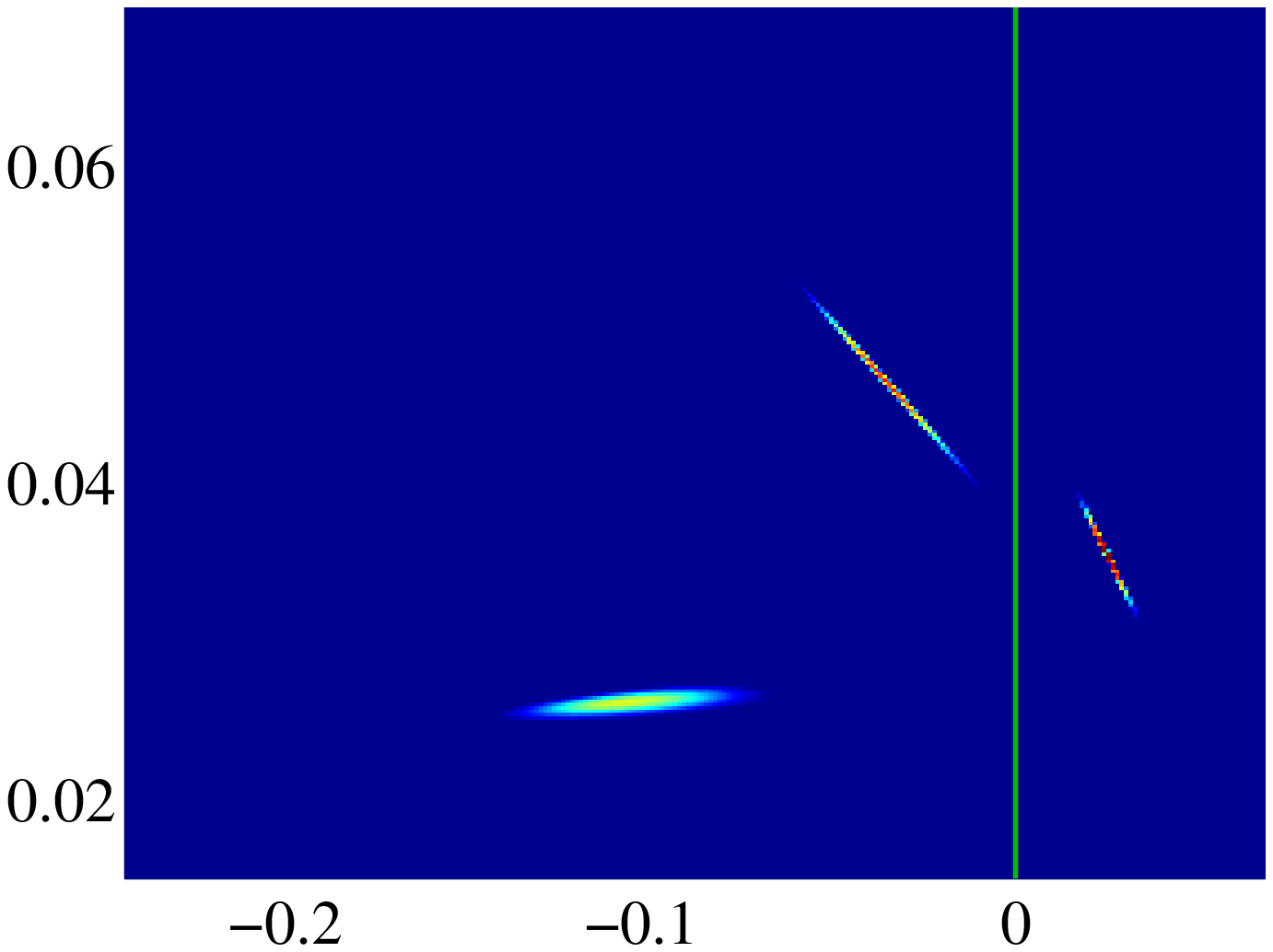}}
\put(7,4.7){\includegraphics[height=5.1cm]{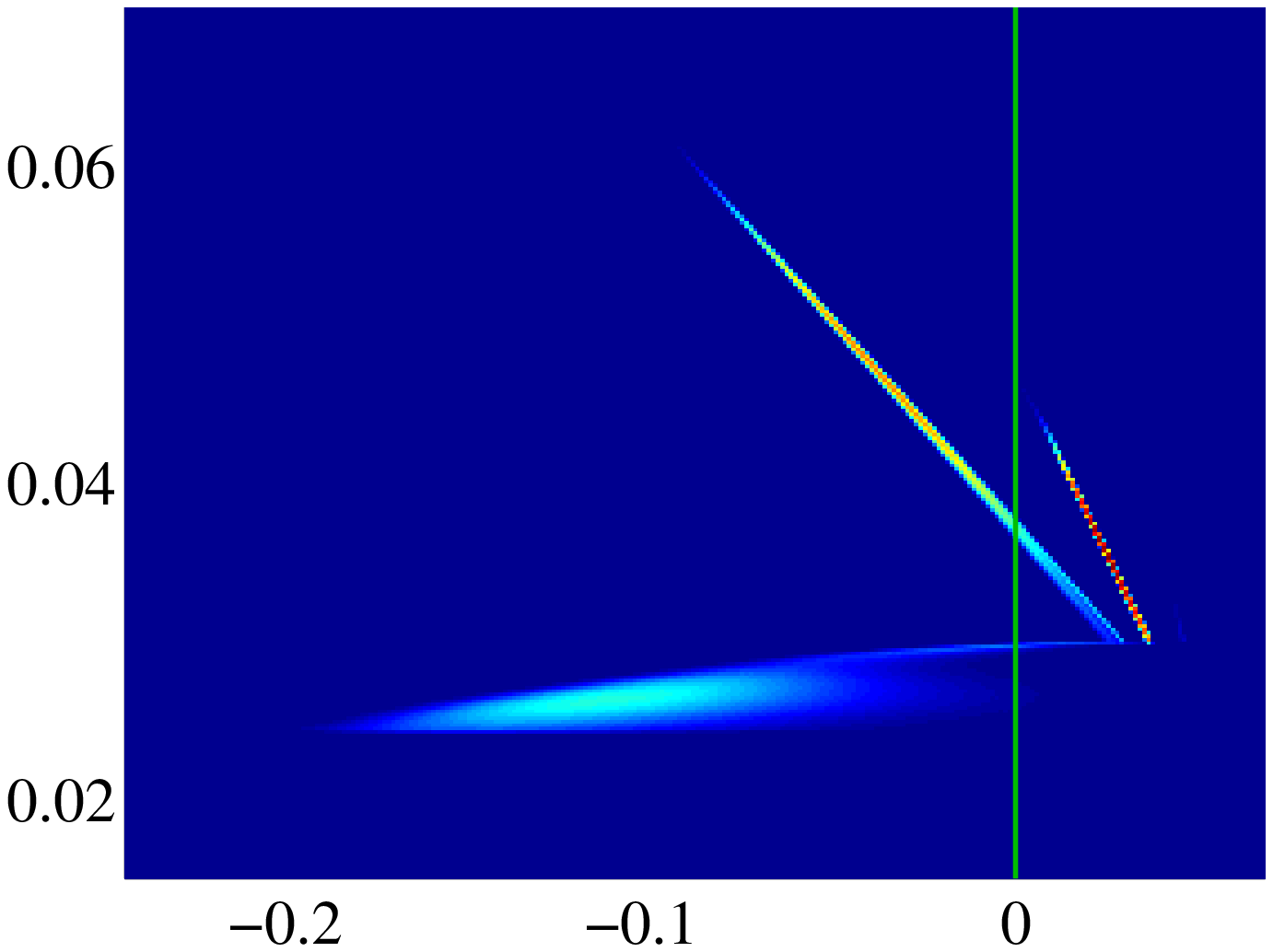}}
\put(2.4,0){\includegraphics[height=4.5cm]{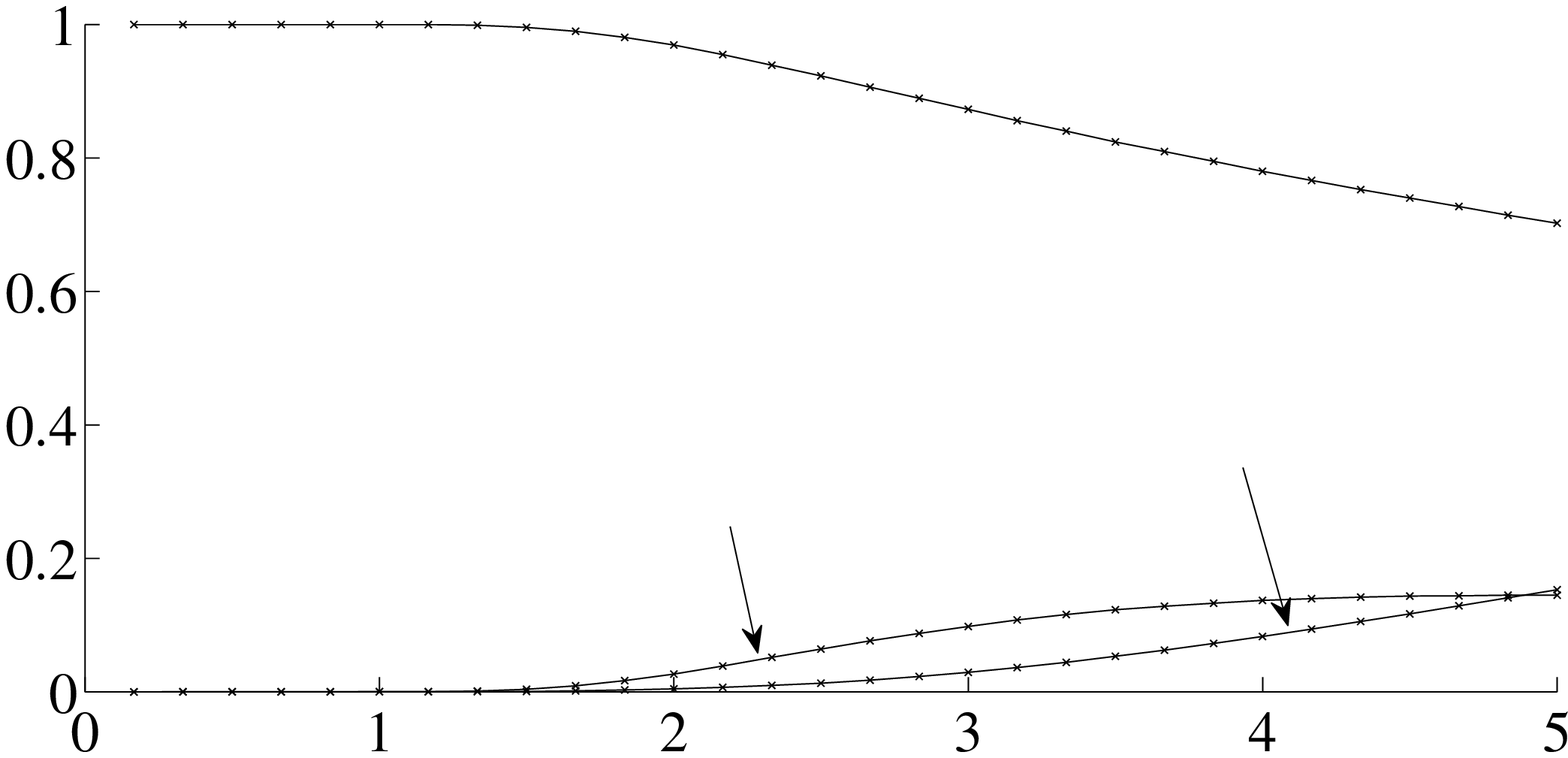}}
\put(3.42,10){\small $\alpha = 1$}
\put(10.42,10){\small $\alpha = 3$}
\put(3.76,4.7){\small $x$}
\put(.3,8.2){\small $y$}
\put(10.76,4.7){\small $x$}
\put(7.3,8.2){\small $y$}
\put(7.2,0){\small $\alpha$}
\put(8,3.44){\small $\sigma_3$}
\put(6.66,1.64){\small $\sigma_2$}
\put(8.78,2){\small $1 - \sigma_2 - \sigma_3$}
\put(.5,9.9){\large \sf \bfseries A}
\put(7.5,9.9){\large \sf \bfseries B}
\put(2.4,4.2){\large \sf \bfseries C}
\end{picture}
\caption{
Panel A shows the invariant density of $N_2$
with the same parameter values as Fig.~\ref{fig:noisyBifDiag} and $\mu = 0.03$.
Panel B shows the invariant density with $\alpha = 3$,
and panel C is a plot of the fractions $\sigma_j$ (\ref{eq:sigma}).
\label{fig:invDensityB_03}
}
\end{center}
\end{figure}

\begin{figure}[b!]
\begin{center}
\setlength{\unitlength}{1cm}
\begin{picture}(13.8,10.2)
\put(0,4.7){\includegraphics[height=5.1cm]{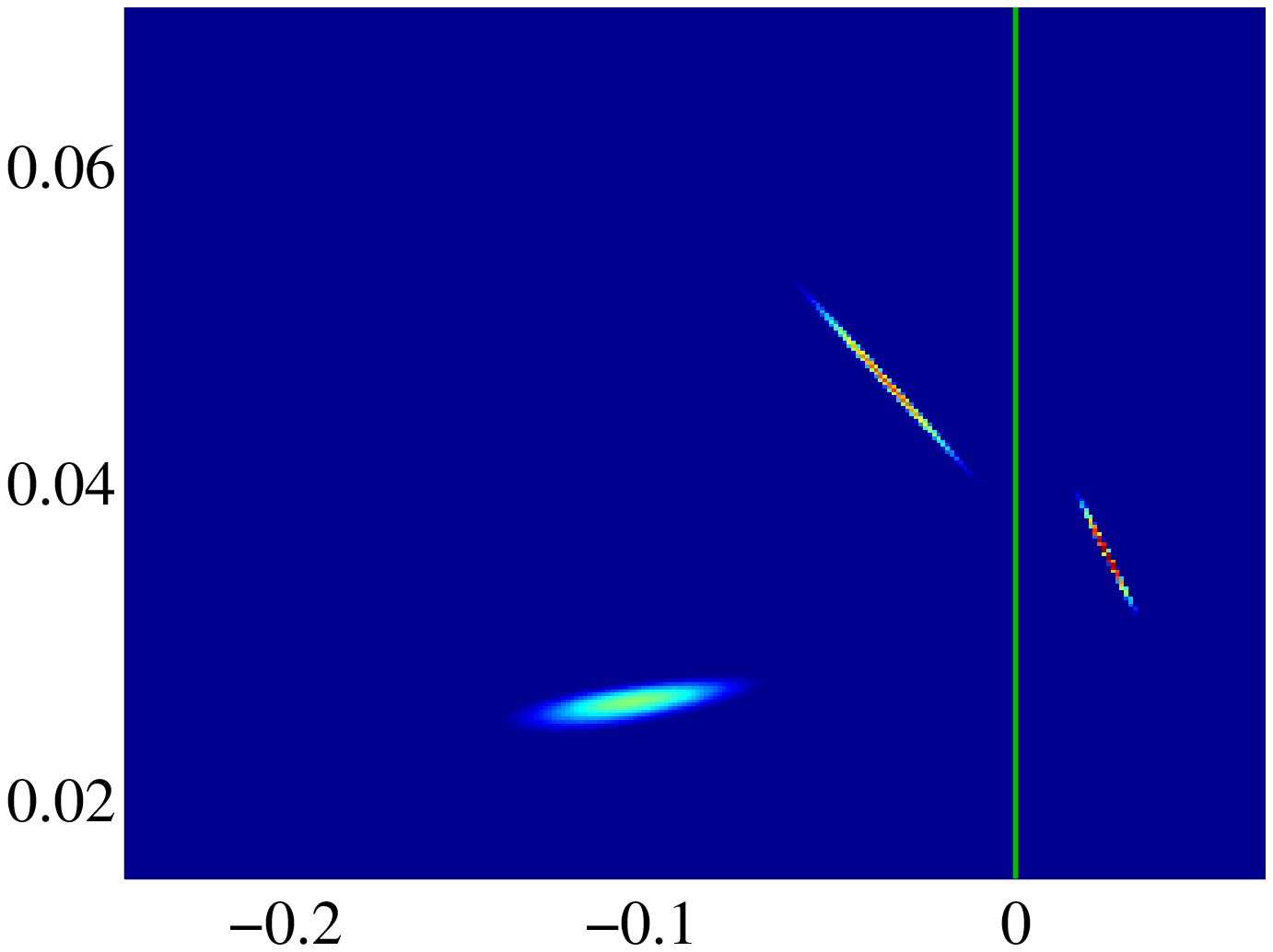}}
\put(7,4.7){\includegraphics[height=5.1cm]{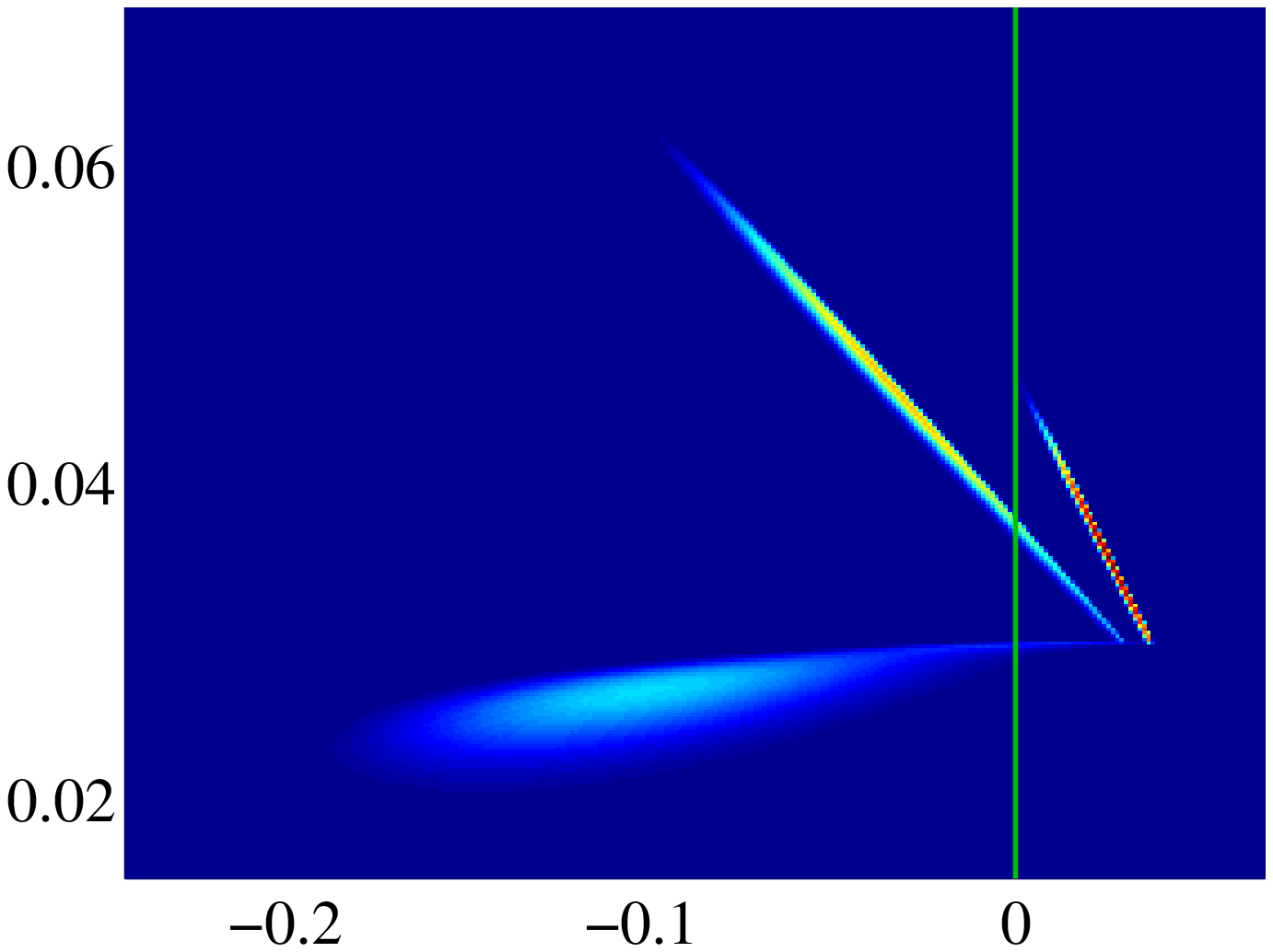}}
\put(2.4,0){\includegraphics[height=4.5cm]{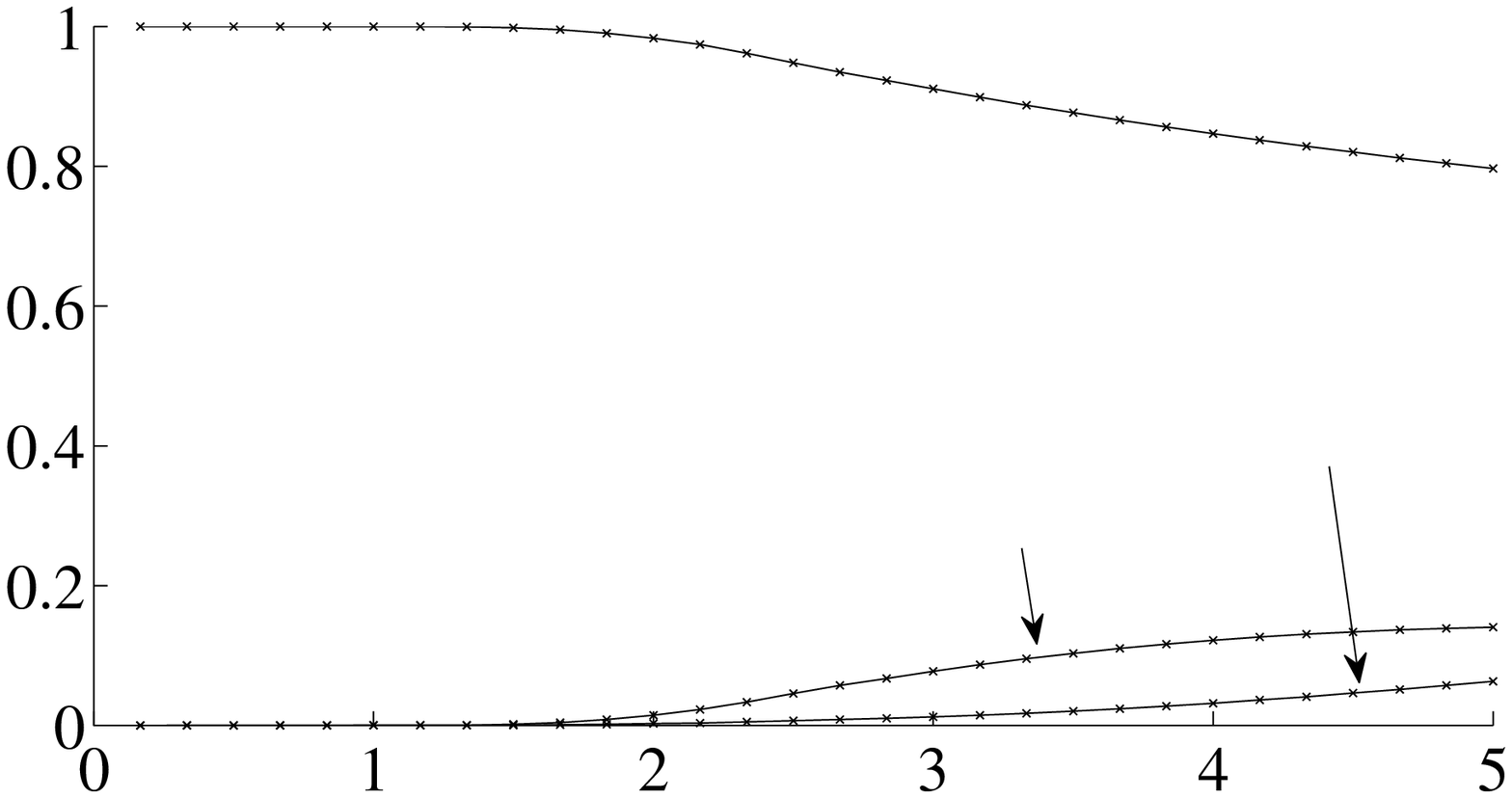}}
\put(3.42,10){\small $\alpha = 1$}
\put(10.42,10){\small $\alpha = 3$}
\put(3.76,4.7){\small $x$}
\put(.3,8.2){\small $y$}
\put(10.76,4.7){\small $x$}
\put(7.3,8.2){\small $y$}
\put(7.2,0){\small $\alpha$}
\put(7.3,3.55){\small $\sigma_3$}
\put(8,1.7){\small $\sigma_2$}
\put(8.93,2.03){\small $1 - \sigma_2 - \sigma_3$}
\put(.5,9.9){\large \sf \bfseries A}
\put(7.5,9.9){\large \sf \bfseries B}
\put(2.4,4.2){\large \sf \bfseries C}
\end{picture}
\caption{
Panel A shows the invariant density of $N_3$
with the same parameter values as Fig.~\ref{fig:noisyBifDiag} and $\mu = 0.03$.
Panel B shows the invariant density with $\alpha = 3$,
and panel C is a plot of the fractions $\sigma_j$ (\ref{eq:sigma}).
\label{fig:invDensityC_03}
}
\end{center}
\end{figure}

Figs.~\ref{fig:invDensityA_03}-\ref{fig:invDensity_m002} show two-dimensional invariant densities of
$N_1$, $N_2$ and $N_3$.
By assuming ergodicity, invariant densities were computed on a $256 \times 256$ grid of $x$ and $y$ values
from $10^8$ consecutive iterates of a single orbit with transient points omitted.

Let us first explain panel C of Figs.~\ref{fig:invDensityA_03}-\ref{fig:invDensityC_03}.
Given a sample orbit 
of $N_1$, $N_2$ or $N_3$,
for each point with $x_i > 0$
we let $\cI_i$ be the smallest positive integer for which $x_{i + \cI_i} > 0$, as in \cite{SiHo13}.
$\cI_i$ represents the number of iterations required for a return to $x > 0$ from the point $(x_i,y_i)$.
Numerically we can compute a large number of values of $\cI_i$ (say $M$ of them).
Then for each $j \in \mathbb{Z}^+$, we let $\sigma_j$ denote the fraction of the $\cI_i$ that are equal to $j$, i.e.
\begin{equation}
\sigma_j = \frac{1}{M} \sum_{i {\rm \,with\,} x_i > 0} \chi_{j - \cI_i} \;,
\label{eq:sigma}
\end{equation}
where $\chi_k = 1$ if $k = 0$, and $\chi_k = 0$ otherwise.
Figs.~\ref{fig:invDensityA_03}-\ref{fig:invDensityC_03}
correspond to $\mu = 0.03$ for which there is an underlying attracting $3$-cycle.
Therefore with small noise, $\sigma_3 \approx 1$, and for all $j \ne 3$, $\sigma_j \approx 0$.

Fig.~\ref{fig:invDensityA_03} corresponds to the map $N_1$.
In panel A, the size of the noise is relatively small,
so iterates of $N_1$ follow close to the $3$-cycle.
The invariant density is well-approximated by a scaled sum of three Gaussian densities
centred at each point of the $3$-cycle.
About the point with $x \approx -0.1$, the density is stretched substantially
more in $x$-direction than in the $y$-direction.
This is because points with $x \approx -0.1$
have likely just undergone an iteration under the right half-map of $N_1$
which is stochastic in the $x$-component but not the $y$-component.
The stretching around other iterates is then a consequence of iterating under $N_1$ with $x < 0$.

With larger values of $\alpha$, it is relatively common for
the orbit to return to $x > 0$ in a number of iterations other than three, Fig.~\ref{fig:invDensityA_03}-C.
For this reason the invariant density displays additional characteristics.
For instance with $\alpha = 3$, the orbit returns to $x > 0$ in two iterations almost $10\%$ of the time.
Consequently, a substantial part of the invariant density
centred roughly about the point of the $3$-cycle with $x \approx -0.04$,
lies in $x > 0$, Fig.~\ref{fig:invDensityA_03}-B.
The invariant density in panel B also has a small component with $x \approx 0.05$.
This is due to points of the orbit with small values of $x > 0$
mapping under the left half-map of $N_1$ due to the noise
(i.e.~returning to $x > 0$ in only one iteration).

Fig.~\ref{fig:invDensityB_03} illustrates $N_2$ using the same parameter values.
Again with small noise the invariant density is roughly Gaussian about each point of the $3$-cycle,
whereas for relatively large noise iterates often cross into $x > 0$ prematurely
causing the invariant density to take an irregular shape.
When $\alpha = 3$, points of the orbit that do not return to $x > 0$ in three iterations,
usually return to $x > 0$ in two iterations.

Fig.~\ref{fig:invDensityC_03} corresponds to $N_3$ and is similar to the previous figure.
This indicates that the correlation time $\nu$ has little effect on these pictures,
although the invariant density has a slightly different shape when $\alpha = 3$.

\subsection{Invariant densities near the grazing bifurcation}

\begin{figure}[b!]
\begin{center}
\setlength{\unitlength}{1cm}
\begin{picture}(13.8,10.4)
\put(0,5.3){\includegraphics[height=5.1cm]{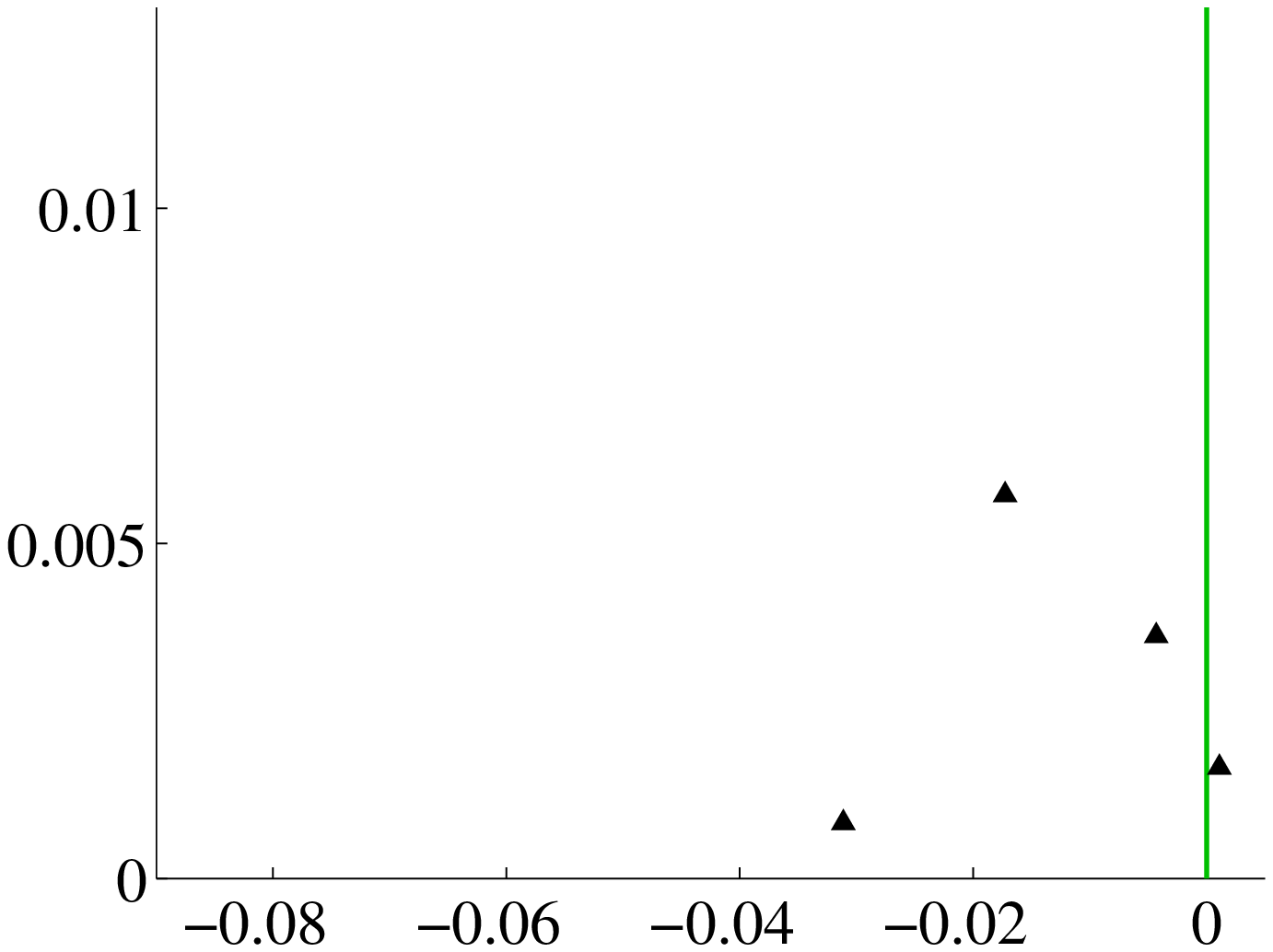}}
\put(7,5.3){\includegraphics[height=5.1cm]{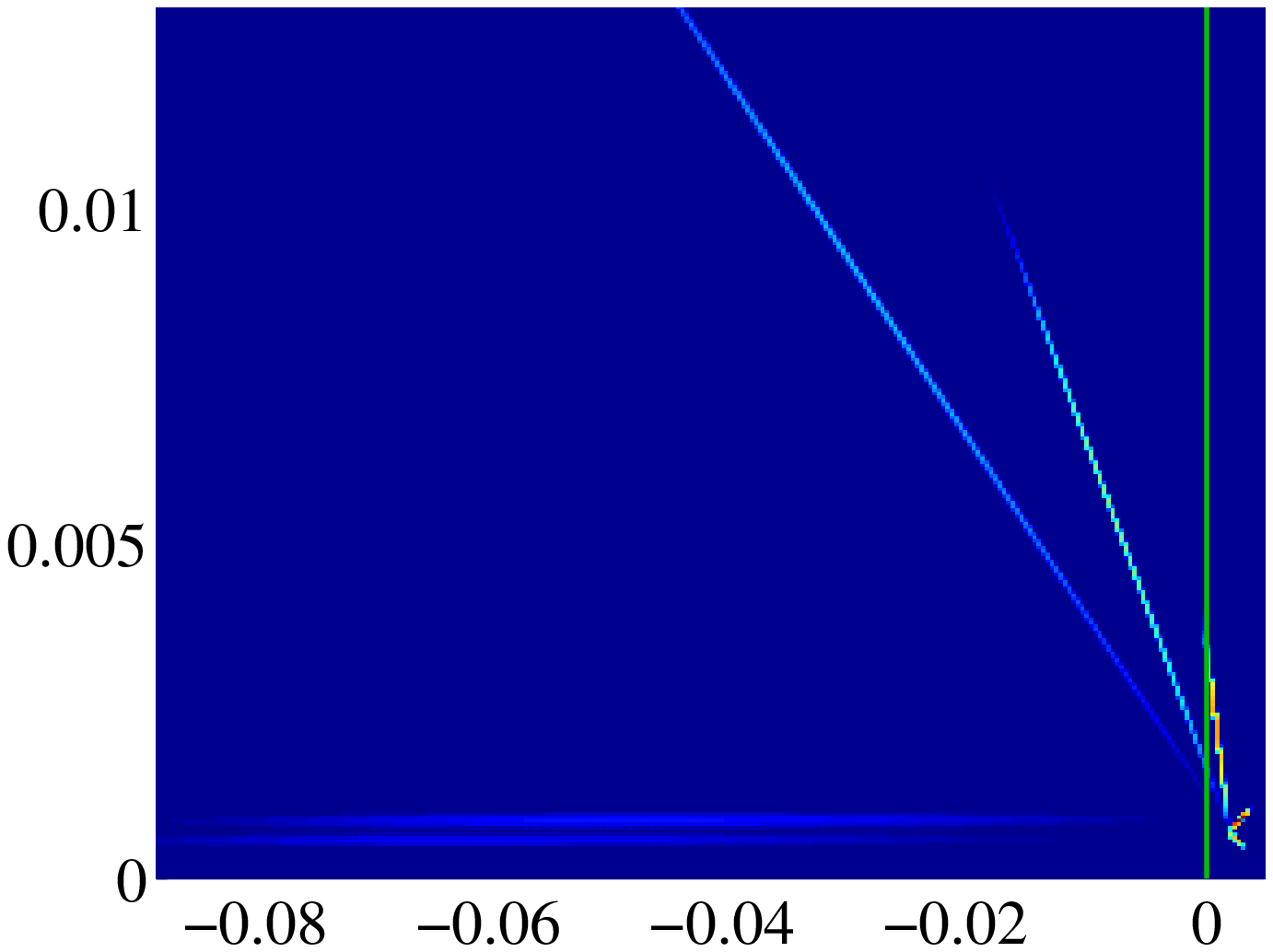}}
\put(0,0){\includegraphics[height=5.1cm]{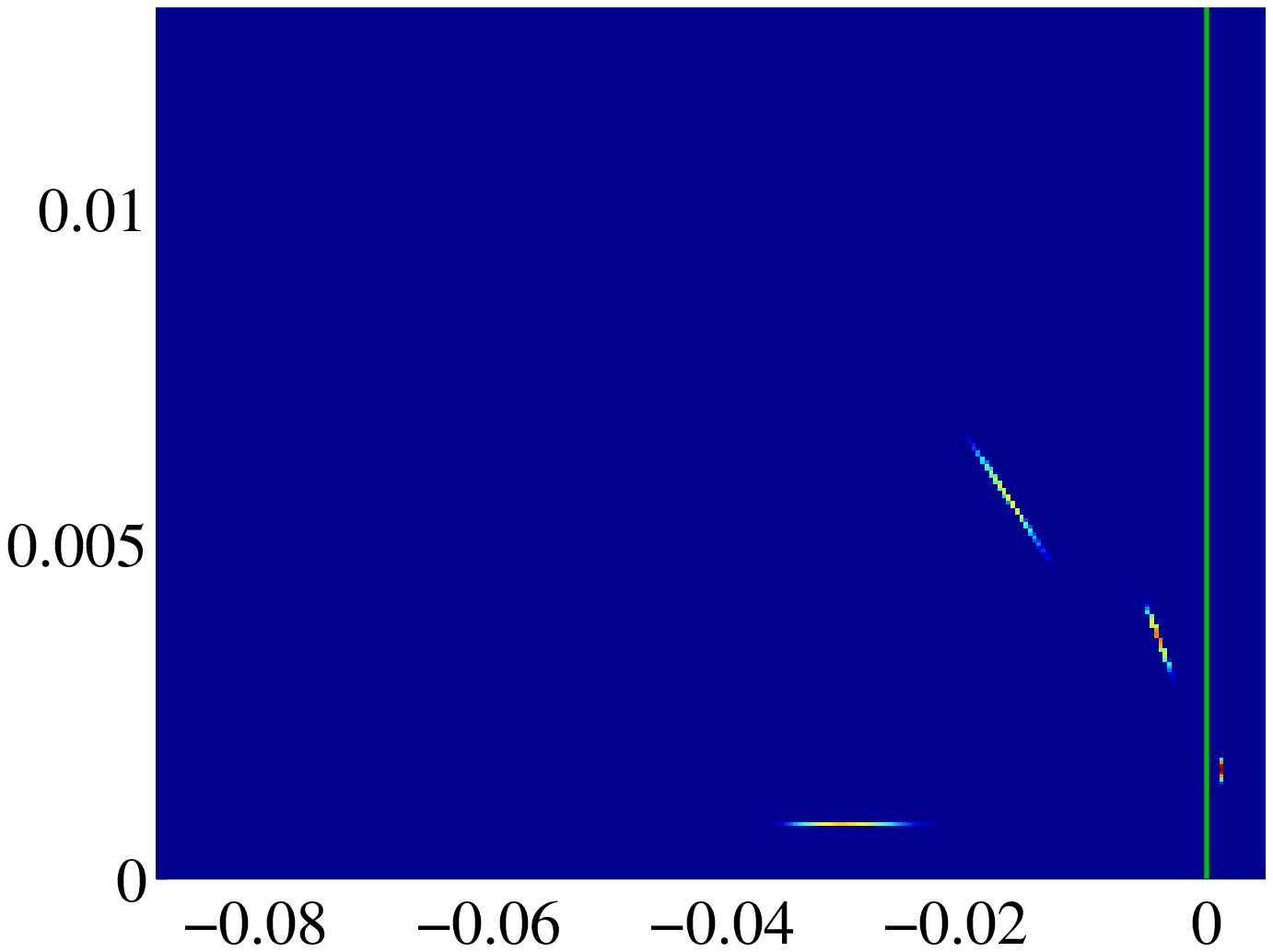}}
\put(7,0){\includegraphics[height=5.1cm]{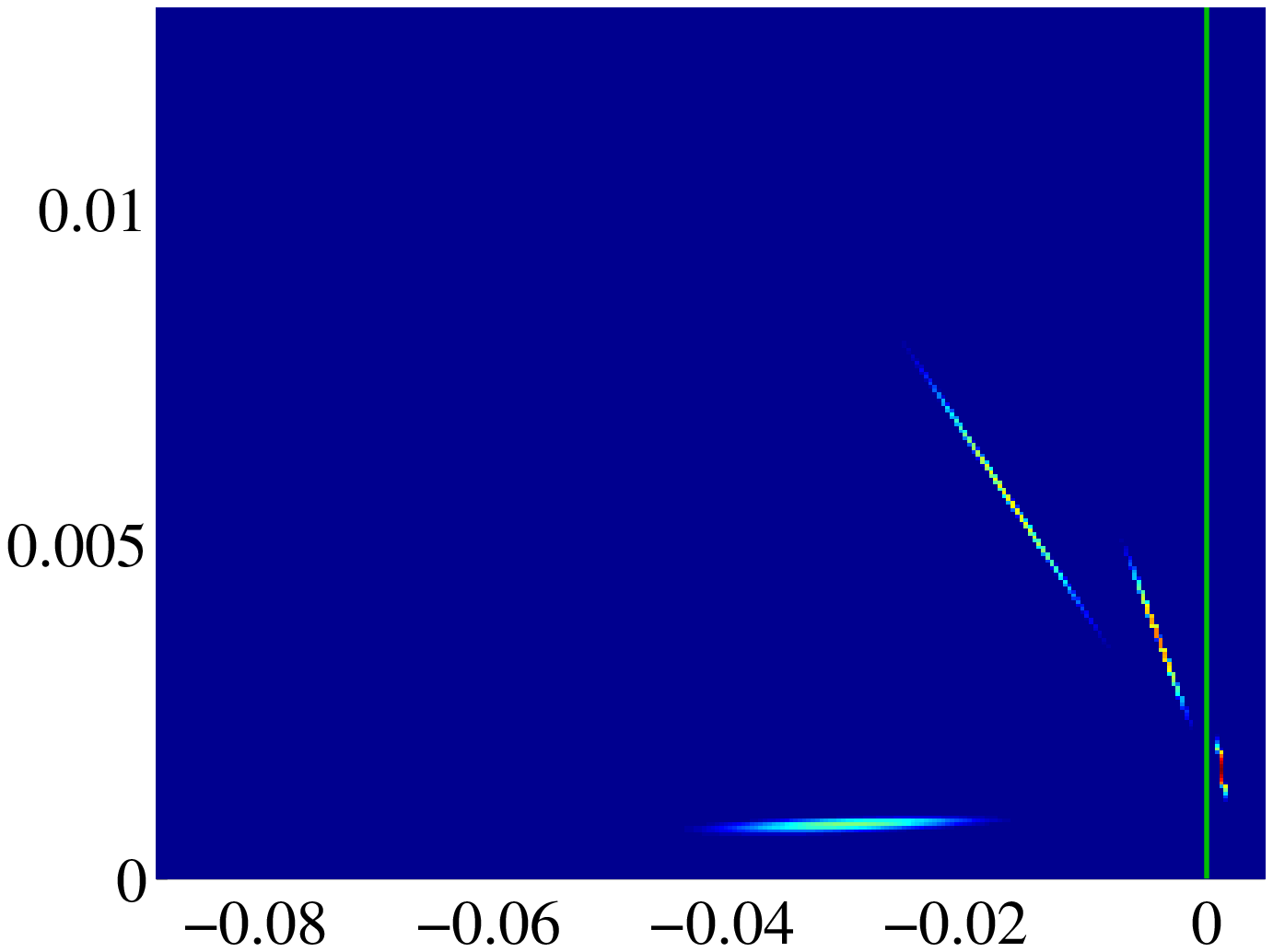}}
\put(3.76,5.3){\small $x$}
\put(.3,8.8){\small $y$}
\put(10.76,5.3){\small $x$}
\put(7.3,8.8){\small $y$}
\put(3.76,0){\small $x$}
\put(.3,3.5){\small $y$}
\put(10.76,0){\small $x$}
\put(7.3,3.5){\small $y$}
\put(.6,10.3){\large \sf \bfseries A}
\put(7.6,10.3){\large \sf \bfseries B}
\put(.6,5){\large \sf \bfseries C}
\put(7.6,5){\large \sf \bfseries D}
\end{picture}
\caption{
Panel A shows the attracting $4$-cycle of (\ref{eq:N})
with the same parameter values as Fig.~\ref{fig:detBifDiag} and $\mu = 0.001$.
Panels B, C and D show the invariant densities of $N_1$, $N_2$ and $N_3$, respectively,
using the same parameter values as Fig.~\ref{fig:noisyBifDiag} and $\mu = 0.001$.
\label{fig:invDensity_001}
}
\end{center}
\end{figure}

Fig.~\ref{fig:invDensity_001} shows invariant densities of $N_1$, $N_2$ and $N_3$
for parameter values closer to the grazing bifurcation than the previous three figures,
specifically $\mu = 0.001$.
At this value of $\mu$, there is an underlying attracting $4$-cycle, panel A.
The noise amplitudes are given by
(\ref{eq:tildeee1}), (\ref{eq:tildeee2}) and (\ref{eq:tildeee3}) with $\alpha = 1$.
Recall, these amplitudes were chosen such that the size of the noise response
of the three maps is roughly the same for larger values of $\mu$.
Here, however, the size of noise response differs substantially.
In panels C and D, which correspond to $N_2$ and $N_3$ respectively,
the invariant density is approximately a scaled sum of four Gaussians
about each point of the $4$-cycle.
The invariant density in panel D, corresponding to $\nu = 0$,
is noticeably larger than that of panel C.

In panel B, which corresponds to $N_1$,
the noise has a substantial effect because the switching condition of $N_1$ is stochastic,
and many points of $N_1$ fall close to $x=0$. 
The invariant density has a small $C$-shaped component in $x>0$
corresponding to consecutive points of the orbit mapping under the left half-map of $N_1$.
The part of the invariant density for $x < 0$ and $y \approx 0.001$
corresponds to images of points under the right half-map of $N_1$,
and is bimodal because the invariant density has roughly two components in $x > 0$.

\subsection{Invariant densities about coexisting attractors}

\begin{figure}[b!]
\begin{center}
\setlength{\unitlength}{1cm}
\begin{picture}(13.8,10.4)
\put(0,5.3){\includegraphics[height=5.1cm]{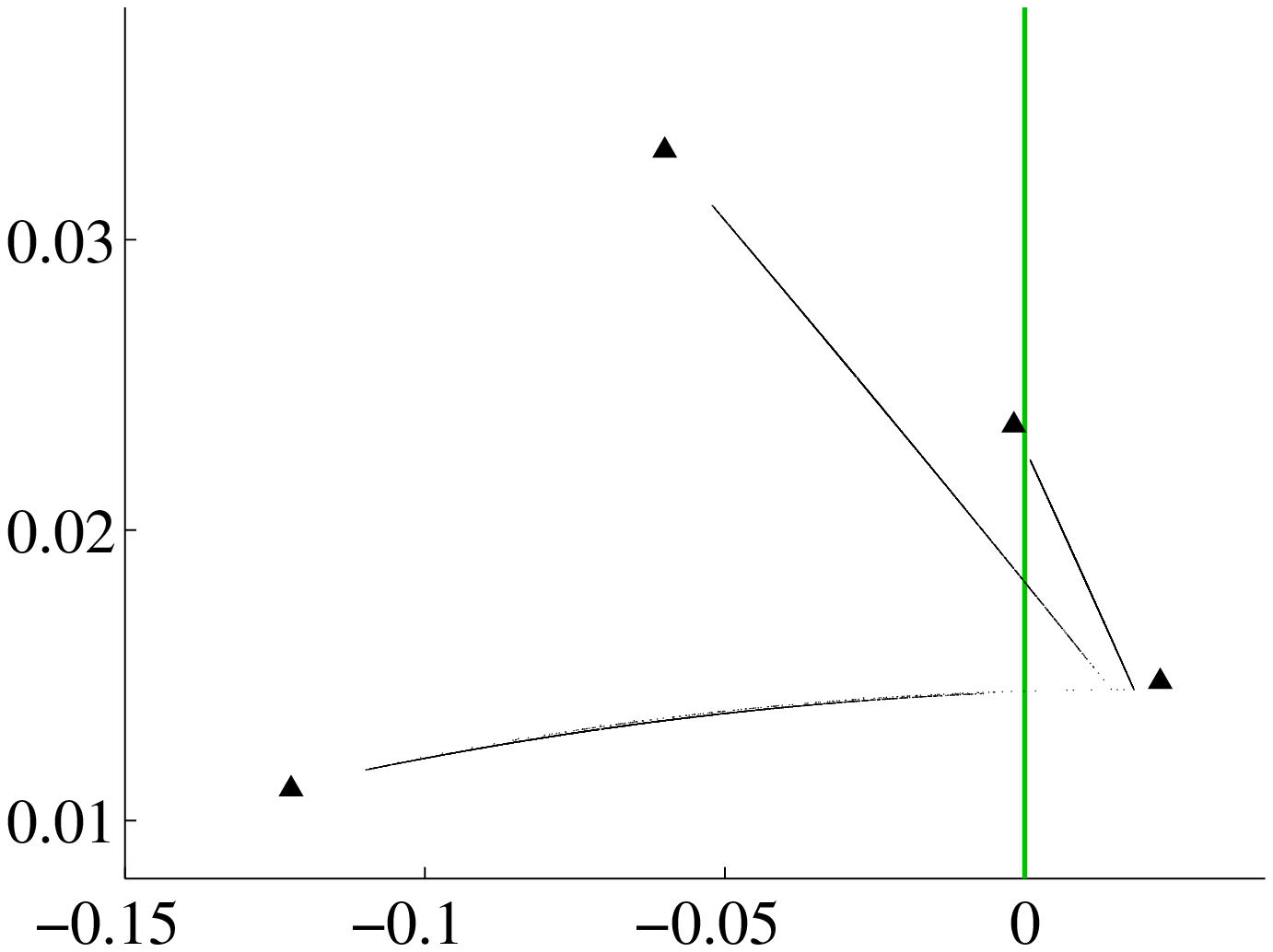}}
\put(7,5.3){\includegraphics[height=5.1cm]{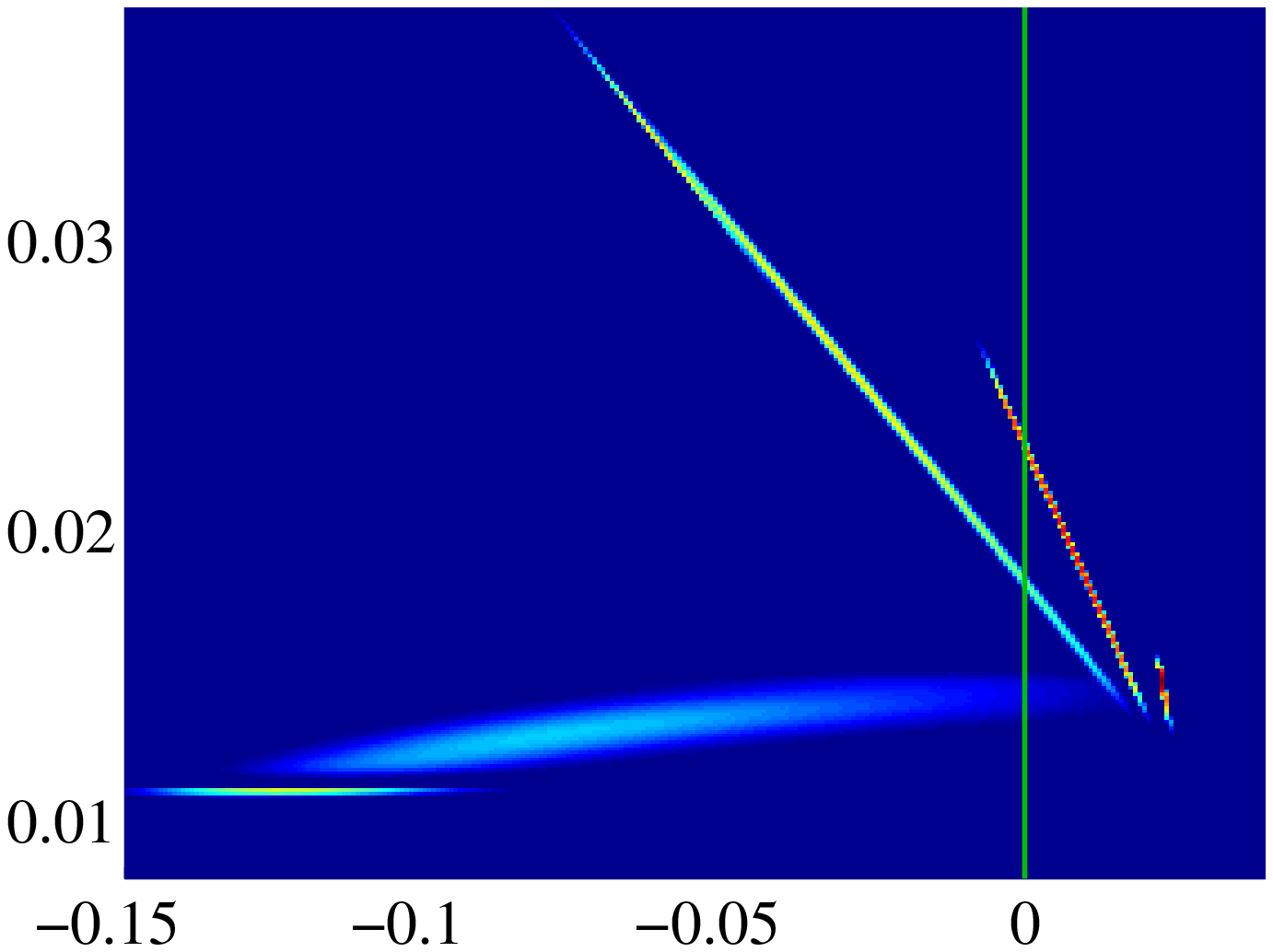}}
\put(0,0){\includegraphics[height=5.1cm]{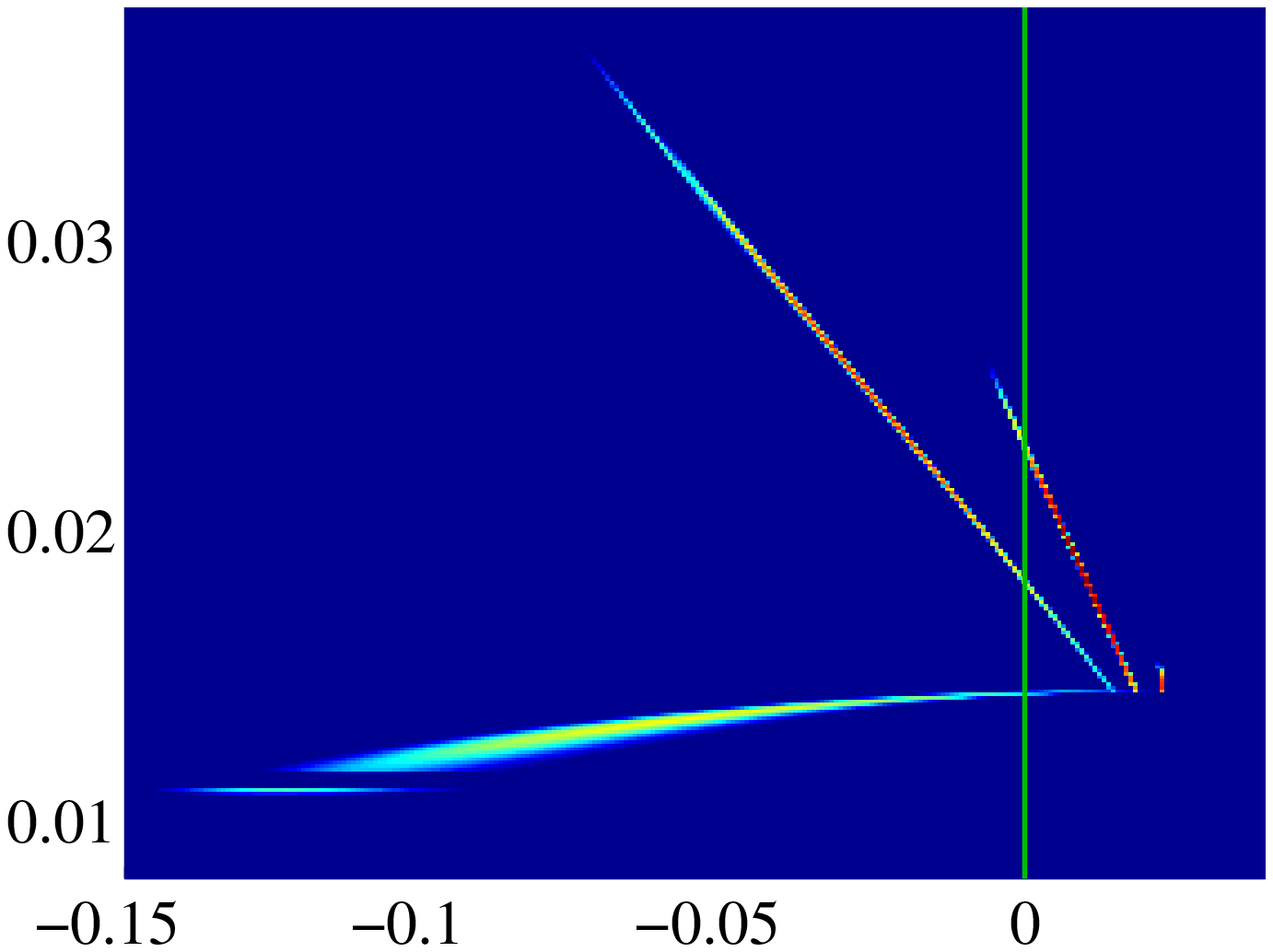}}
\put(7,0){\includegraphics[height=5.1cm]{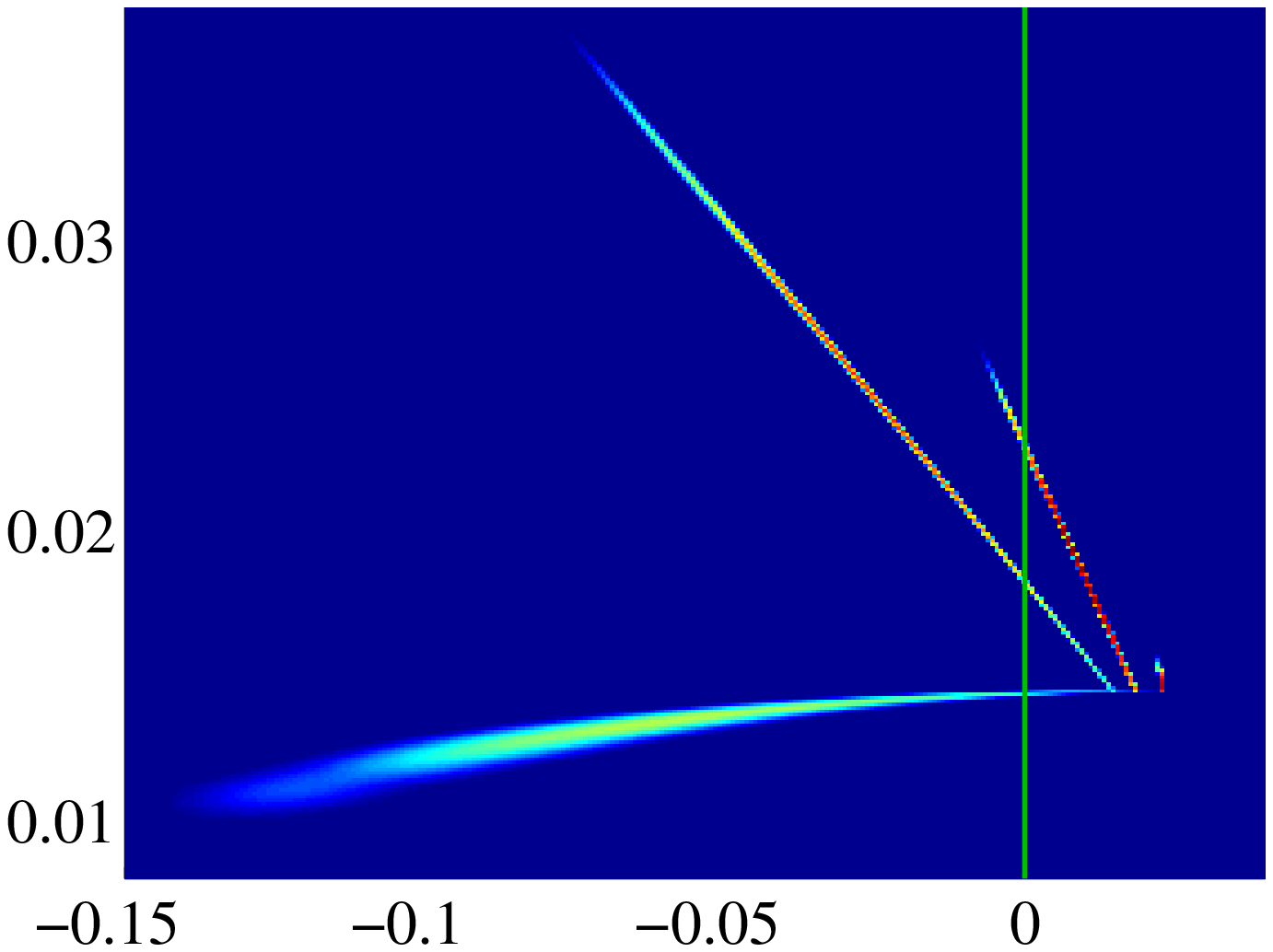}}
\put(3.76,5.3){\small $x$}
\put(.3,8.5){\small $y$}
\put(10.76,5.3){\small $x$}
\put(7.3,8.5){\small $y$}
\put(3.76,0){\small $x$}
\put(.3,3.2){\small $y$}
\put(10.76,0){\small $x$}
\put(7.3,3.2){\small $y$}
\put(.4,10.3){\large \sf \bfseries A}
\put(7.4,10.3){\large \sf \bfseries B}
\put(.4,5){\large \sf \bfseries C}
\put(7.4,5){\large \sf \bfseries D}
\end{picture}
\caption{
Panel A shows an attracting $4$-cycle
and a numerically computed attractor of (\ref{eq:N})
with the same parameter values as Fig.~\ref{fig:detBifDiag} and $\mu = 0.0145$.
Panels B,C and D show the invariant densities of $N_1$, $N_2$ and $N_3$, respectively,
with the same parameter values as Fig.~\ref{fig:noisyBifDiag} and $\mu = 0.0145$.
\label{fig:invDensity_0145}
}
\end{center}
\end{figure}

\begin{figure}[t!]
\begin{center}
\setlength{\unitlength}{1cm}
\begin{picture}(13.8,5.1)
\put(0,0){\includegraphics[height=5.1cm]{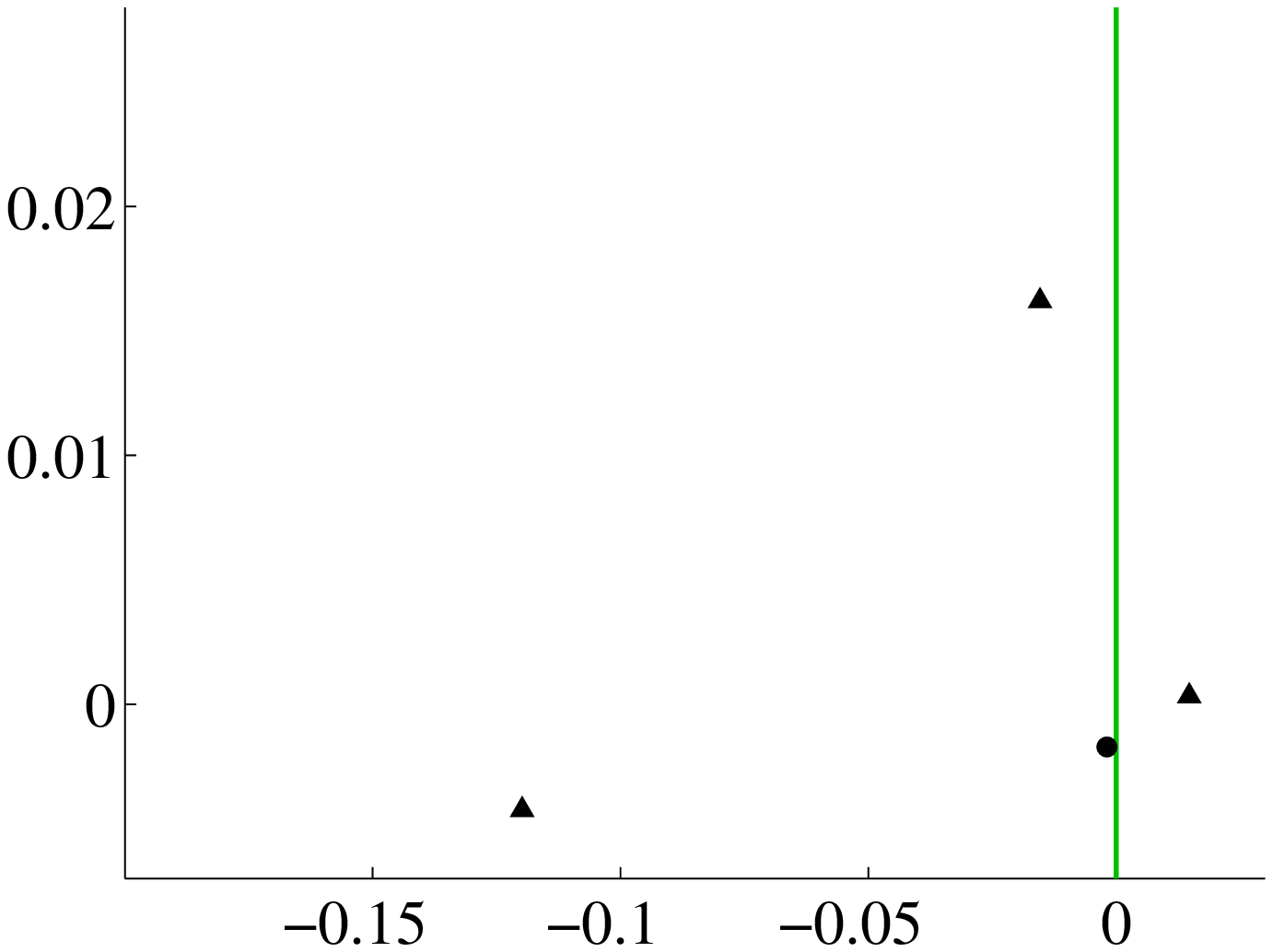}}
\put(7,0){\includegraphics[height=5.1cm]{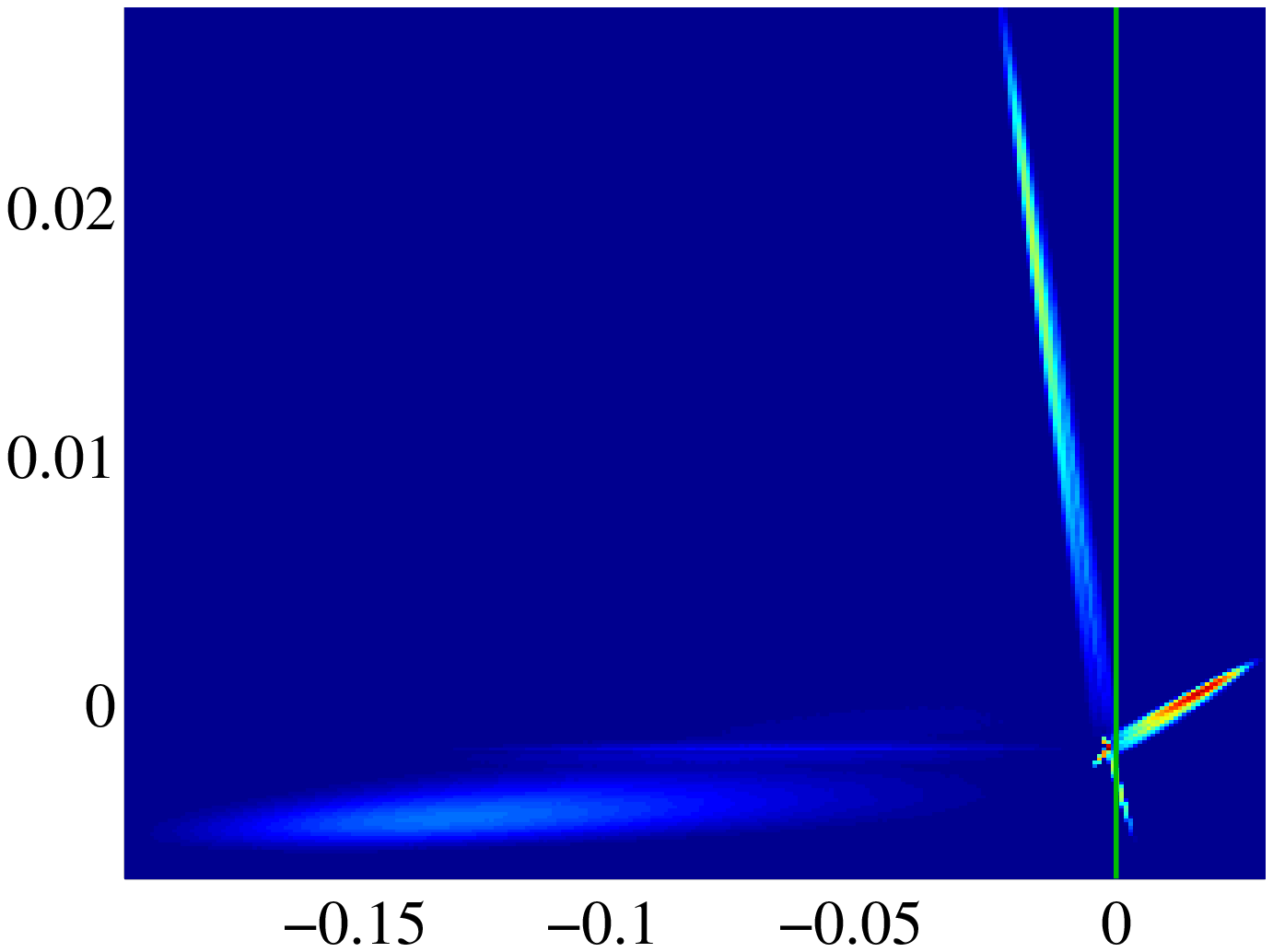}}
\put(3.76,0){\small $x$}
\put(.3,3.4){\small $y$}
\put(10.76,0){\small $x$}
\put(7.3,3.4){\small $y$}
\put(.4,5){\large \sf \bfseries A}
\put(7.4,5){\large \sf \bfseries B}
\end{picture}
\caption{
Panel A shows the attracting fixed point $(x^L,y^L)$ (\ref{eq:xLyL})
and the attracting $3$-cycle of $N$
with the same parameter values as Fig.~\ref{fig:detBifDiag}, except $k_{\rm osc} = 5$ and $\mu = -0.002$.
Here $\tau \approx 0.0927$, $\delta \approx 0.1518$ and $\chi = 1$.
Panel B shows the invariant density of the corresponding map $N_1$, with $\nu = 0.5$ and $\ee = 0.0005$.
\label{fig:invDensity_m002}
}
\end{center}
\end{figure}

With $\mu = 0.0145$ in Fig.~\ref{fig:detBifDiag},
there is an attracting $4$-cycle and an apparently chaotic attractor.
These are shown in Fig.~\ref{fig:invDensity_0145}-A.
As with smooth maps \cite{KnWe89,KrFe99},
in the presence of noise orbits commonly dwell near the attractors for relatively long periods of time,
and switch between attractors quickly.
Invariant densities of $N_1$, $N_2$ and $N_3$ are shown in panels B, C and D. 
In each case the bulk of the density is centred about the two underlying attractors.
With white noise (panel D) there is no gap in the invariant density
around $(x,y) \approx (-0.12,0.012)$ due to randomness in both the $x$ and $y$-components of $N_3$.

Lastly, Fig.~\ref{fig:invDensity_m002} illustrates stochastic dynamics with $\mu < 0$.
This figure corresponds to $k_{\rm osc} = 5$ (different to Fig.~\ref{fig:detBifDiag})
and $\mu = -0.002$ at which $N$ has an attracting $3$-cycle
as well as the attracting fixed point $(x^L,y^L)$ (\ref{eq:xLyL}).
These are shown in panel A.
The dynamics of $N_2$ and $N_3$ for $x < 0$ is deterministic,
hence $(x^L,y^L)$ is a fixed point of these maps.
Given an initial point $(x_0,y_0)$ near the $3$-cycle,
sample orbits of $N_2$ and $N_3$ eventually reach $(x^L,y^L)$\removableFootnote{
Escape times should be geometrically distributed,
and the mean escape time should decrease exponentially with the noise amplitude,
see for instance \cite{Be89,HiMe13}.
}.
In contrast, $N_1$ has an invariant density
concentrated about the two attractors, panel B.
The part of the density with $x \approx 0$ and $y < y_L$
corresponds to points of the orbit repeatedly following the left half-map of $N_1$
(with $x_i < 0$ and $x_i + \frac{\xi_i}{a_{12}^2 c^2} < 0$).

\section{Conclusions}
\label{sec:conc}
\setcounter{equation}{0}

This paper concerns grazing bifurcations
for which the associated dynamics is described by the Nordmark map (\ref{eq:N}).
The potential influence of randomness and uncertainties on the dynamics of (\ref{eq:N})
was described in \cite{SiHo13} by studying (\ref{eq:N}) in the presence of additive white Gaussian noise.
Such a noise formulation is suitable if the nature of the randomness in the ODE system
is practically independent to the state of the system,
such as if there is a random forcing term.

In this paper we considered the alternate scenario
that randomness is present in the switching condition of the ODE system,
and in $f_R$ -- the part of the vector field opposite to the tangential intersection of the grazing periodic orbit.
These cases are especially relevant for vibro-impacting systems
for which impact events represent the primary source of uncertainty.
We derived three different stochastic versions of (\ref{eq:N}).
These are the maps $N_1$ (\ref{eq:Na}), which corresponds to a noisy switching condition in the ODE system,
$N_2$ (\ref{eq:Nb}), which corresponds to noise in $f_R$ with a large correlation time,
and $N_3$ (\ref{eq:Nc}), which corresponds to white noise in $f_R$.
In each case the noise is nonlinear and non-additive.
This indicates that some diligence should be taken when formulating stochastic return maps
for grazing bifurcations of piecewise-smooth systems.

The stochastic dynamics of $N_1$, $N_2$ and $N_3$
differs in many ways to that of (\ref{eq:N}) with additive noise, described in \cite{SiHo13}.
For $N_1$, $N_2$ and $N_3$, dynamics prior to the grazing bifurcation is deterministic,
and beyond the grazing bifurcation two-dimensional invariant densities
are often skewed dramatically so that they appear almost one-dimensional.

Near the grazing bifurcation, $N_1$ exhibits a large noise response relative to $N_2$ and $N_3$.
This suggests that if experimental data of a physical system shows relatively high variability near a grazing bifurcation,
then it is likely to be most appropriate to include randomness in the switching condition of a mathematical model.
Invariant densities of $N_1$ near grazing are highly irregular due to the randomness in the switching condition.
For $N_2$ and $N_3$, the size of the noise response increases, for most part,
with the distance (in parameter space) beyond the grazing bifurcation.
The maps $N_2$ and $N_3$ exhibit qualitatively similar invariant densities,
which implies that the correlation time has little effect.
Indeed the correlation time only influences the short-time dynamics of (\ref{eq:fbc}) with $u > 0$,
and the precise nature of these dynamics has a negligible effect on the invariant densities, which relate to long-time dynamics.

\end{document}